\documentclass[a4paper,12pt]{article}
\usepackage[centertags]{amsmath}
\usepackage[french,english]{babel}
\usepackage{amsfonts}
\usepackage{amsmath,amssymb,graphics,amscd}
\usepackage{amsthm}
\input pictexwd

\newcommand{\semi}{{\,\rule[.1pt]{.4pt}{5.3pt}\hskip-1.9pt\times}}
\newcommand{\la}{\langle}
\newcommand{\ra}{\rangle}
\newcommand{\ce}{\widetilde{e}}
\newcommand{\al}{\alpha}
\newcommand{\be}{\beta}
\newcommand{\alvee}{{\al^\vee}}

\newcommand{\Om}{{\Omega}}
\newcommand{\g}{\gamma}
\newcommand{\G}{\Gamma}
\newcommand{\BSlam}{\hat\Sigma(\g_\lam )}
\newcommand{\de}{\delta}

\newcommand{\De}{\Delta}
\newcommand{\M}{{\mathcal M}}

\newcommand{\A}{\mathbb A}

\newcommand{\set}[1]{\left\{#1\right\}}

\newcommand{\avee}{\alpha^\vee}

\newcommand{\Xvee}{{X^\vee}}
\newcommand{\Xveep}{{X^\vee_+}}
\newcommand{\Phivee}{{\Phi^\vee}}
\newcommand{\Phiveep}{{\Phi^\vee_+}}

\newcommand{\Rvee}{{R^\vee}}
\newcommand{\Rveep}{{R^\vee_+}}
\newcommand{\Delvee}{\Delta^\vee}
\newcommand{\Waff}{{W^{{\mathfrak a}}}}
\newcommand{\Saff}{{S^{{\mathfrak a}}}}
\newcommand{\Aaff}{{A^{{\mathfrak{a}}}}}

\newcommand{\lam}{{\lambda}}

\newcommand{\wte}{{\widetilde{e}}}
\newcommand{\rinf}{r_{-\infty}}

\newcommand{\cu}{\mathcal{U}}
\newcommand{\co}{\mathcal{O}}
\newcommand{\ck}{\mathcal{K}}
\newcommand{\cg}{\mathcal{G}}
\newcommand{\cl}{\mathcal{L}}

\newcommand{\cj}{\mathcal{J}}
\newcommand{\cjaff}{\mathcal{J}^{\mathfrak a}}

\newcommand{\cb}{\mathcal{B}}

\newcommand{\cp}{\mathcal{P}}
\newcommand{\cq}{\mathcal{Q}}
\newcommand{\ca}{\mathcal{A}}
\newcommand{\tta}{\mathcal{A}}

\newcommand{\ch}{{\mathtt{H}}}
\newcommand{\cha}{{\mathtt{H}^{\mathfrak a}}}

\newcommand{\cn}{\mathcal{N}}

\newcommand{\ct}{\mathcal{T}}
\newcommand{\sh}{\mathcal{H}}

\newcommand{\chara}{{\mathop{\rm Char}\,}}

\newcommand{\Lie}{{\mathop{\hbox{\rm Lie}}\,}}

\newcommand{\Mor}{{{\hbox{\rm Mor}}\,}}
\newcommand{\Aut}{{{\hbox{\rm Aut}}\,}}

\newcommand{\sgn}{{\mathop{\rm sgn}}}

\newcommand{\br}{{\mathbb R}}
\newcommand{\bz}{{\mathbb Z}}

\newcommand{\bc}{{\mathbb C}}

\newcommand{\Lg}{{\mathfrak{g}}}

\newcommand{\Lc}{{\mathfrak{C}}}
\newcommand{\Ls}{{\mathfrak{s}}}
\newcommand{\Lcinf}{{\mathfrak{C}_{-\infty}}}
\newtheorem{thm}{Theorem}
\newtheorem{dfn}{Definition}
\newtheorem{defiprop}{Definition-Proposition}
\newtheorem{exam}{Example}
\newtheorem{rem}{Remark}
\newtheorem{rems}[rem]{Remarks}
\newtheorem{prop}{Proposition}
\newtheorem{lem}{Lemma}
\newtheorem{coro}{Corollary}
\def\mapright#1{\smash{\mathop{\longrightarrow}\limits^{#1}}}

\def\proof{{\it Proof. }}
\def\endpf{\hfill$\bullet$\medskip}

\begin{document}
\title{LS-Galleries, the path model and MV-cycles}
\author{S. Gaussent$^*$ and P. Littelmann\footnote{This research has been partially 
supported by the EC TMR network "Algebraic Lie Representations", contract
no. ERB FMRX-CT97-0100. {\sl 2000 Mathematics Subject Classification. 22E46, 14M15.}}}
\maketitle
\bigskip

\begin{abstract}
We give an interpretation of the path model of a representation \cite{Lit1} 
of a complex semisimple algebraic group $G$ in terms of the geometry of its affine Grassmannian. 
In this setting, the paths are replaced by LS--galleries in the affine Coxeter complex associated 
to the Weyl group of $G$. To explain the connection with geometry, consider a
Demazure--Hansen--Bott--Samelson desingularization $\hat\Sigma(\lam)$ of the closure of 
an orbit $G(\bc[[t]]).\lam$ in the affine Grassmannian. The
homology of $\hat\Sigma(\lam)$ has a basis given by Bia{\l}ynicki--Birula cell's, which are indexed by the
$T$--fixed points in $\hat\Sigma(\lam)$. Now the points of $\hat\Sigma(\lam)$ can be 
identified with galleries of a fixed type in the affine Tits building associated to $G$, and the
$T$--fixed points correspond in this language to combinatorial galleries of a fixed type in the affine
Coxeter complex. We determine those galleries such that the associated cell has a
non-empty intersection with $G(\bc[[t]]).\lam$ (identified with an open subset of $\hat\Sigma(\lam)$), 
and we show that the closures of the strata associated to LS-galleries are exactly the 
MV--cycles \cite{MV}, which form a basis of the representation
$V(\lam)$ for the Langland's dual group $G^\vee$.
\end{abstract}
\section*{Introduction}
The aim of the present article is to provide a connection between the path model
\cite{Lit1} for finite dimensional representations of a connected complex semisimple 
algebraic group $G$ and the geometry of the affine Grassmann variety $\cg$. As a 
byproduct of this interpretation, we associate in a canonical way to each path 
(or rather gallery) in a path model of a representation one of the algebraic cycles 
occurring in the work of Mirkovic and Vilonen \cite{MV}. Recall that these cycles form, 
by \cite{MV} and Vasserot \cite{Vass}, a canonical basis of the finite dimensional complex 
representations of the Langland's dual group $G^\vee$. 

The affine Grassmannian, as a set, is defined as
$\cg=G(\ck)/G(\co)$, where $\ck=\bc((t))$ and $\co=\bc[[t]]$. 
Fix a maximal torus $T\subset G$ and let $X=X(T)$ be its character group.
The group $\Xvee = \Mor (\bc^*, T)$ of co-characters
is the character group of the dual maximal torus $T^\vee$ of $G^\vee$.
The dominant co-characters $X^\vee_+$ form an indexing set for two 
classes of objects. On the one hand, the finite dimensional irreducible $G^\vee$-modules $V(\lam)$ 
are classified by their highest weight $\lam\in X^\vee_+$ . On the other hand, the co-characters 
can be viewed as points in $\cg$, and $\cg$ is the disjoint union of the $G(\co)$--orbits 
$\cg_\lam = G(\co)\lam$, $\lam\in X^\vee_+$. The (Zariski) closure 
$X_\lam = \overline{\cg_\lam}$ is a finite dimensional projective variety, we fix a 
Demazure--Hansen--Bott--Samelson desingularization $\hat\Sigma(\lam)\rightarrow X_\lam$.

Another object which comes into the picture is the affine Tits--building $\cjaff$ associated to $G$.
Corresponding to the choice of $T\subset G$, we have 
a fixed apartment $\ca\subset \cjaff$ in the Tits building, which, as a simplicial complex, 
is isomorphic to the affine Coxeter complex associated to the Weyl group of $G$.
The link between the combinatorial part and the geometrical part of our construction
is described by the diagram below:
$$
\begin{array}{ccccccc}
G(\co)\hbox{\rm--orbit\ }\cg_\lam&\hookrightarrow &\hat\Sigma(\lam)\cr
=& &=\cr
\hbox{\rm minimal galleries of fixed type\ } \gamma_\lam&\hookrightarrow&
                                         \hbox{\rm galleries of fixed type} \gamma_\lam\cr 
\hbox{\rm in the affine building\ } \cjaff&
&\hbox{\rm in the affine building\ } \cjaff\cr
\downarrow {{\hbox{\it\ retraction onto the}}\atop{\hbox{\it affine building}}}&
& \downarrow{{\hbox{\it\ $T$--fixed points}}\atop {\hbox{\it in BB-cells}}}\cr
\hbox{\rm positively folded combinatorial}&\hookrightarrow
&\hbox{\rm combinatorial galleries of }\cr
\hbox{\rm galleries of fixed type $\gamma_\lam$ in\ }\ca&
&\hbox{\rm fixed type $\gamma_\lam$ in\ }\ca\cr
\end{array}
$$
which we explain now, starting
with the combinatorial galleries.
Assume for simplicity that $\lam$ is regular. A {\it combinatorial gallery} 
joining the origin and $\lam$ is a 
sequence of alcoves $(\Delta_i)_{i=0,\ldots,r}$ in $\ca$ such that 
$0\in\Delta_0$, $\lam\in\Delta_r$, and two consecutive alcoves have 
at least a codimension 
one face in common. Fix a minimal such gallery $\gamma_\lam$ and 
let $\Gamma(\gamma_\lam)$ 
denote the set of all combinatorial galleries in $\ca$ of the same type as
$\gamma_\lam$ (section~\ref{combinatorial gallery}) 
starting at the origin. We define a dimension function on the set $\Gamma(\gamma_\lam)$
counting the number of certain hyperplanes crossed by a gallery,
see section~\ref{dimensiongallery}.  The dimension of the galleries is bounded 
above, and we call a gallery in $\Gamma(\gamma_\lam)$ an LS--gallery if its 
dimension is equal to the upper bound.
We construct ``folding operators'' $e_\alpha,f_\alpha$ on the set of galleries 
$\Gamma(\gamma_\lam)$ for all simple roots and show by analogy with the path model 
\cite{Lit1}: 
(see Proposition~\ref{characterformula}, Theorem~\ref{crystal})
\par\vskip 8pt\noindent
{\bf Theorem A.} {\it The set of LS--galleries $\Gamma_{LS}(\gamma_\lam)$ in $\Gamma(\gamma_\lam)$ 
is stable under the folding operators. Let $B(\gamma_\lam)$ be the directed colored 
graph having as vertices the set of LS-galleries $\Gamma_{LS}(\gamma_\lam)$, 
and put an arrow $\delta\mapright{\alpha}\delta'$
with color $\alpha$ between two galleries if $f_\alpha(\delta)=\delta'$. Then this graph is
connected, and it is isomorphic to the crystal graph of the irreducible 
representation $V(\lam)$ of $G^\vee$ of highest weight $\lam$.  
Let $e(\delta)$ be the endpoint of a gallery, then
$$
\chara V(\lam)=\sum_{\delta\in \Gamma_{LS}(\gamma_\lam)} \exp(e(\delta)).
$$
}
To  connect the combinatorics  with the geometry, recall that in terms of the 
affine Kac-Moody group associated to the extended Dynkin diagram of $G$, 
the affine Grassmannian is a generalized flag variety and $X_\lam$ is a Schubert 
variety. Fixing a minimal gallery in $\ca$ joining the origin and $\lam$ is 
equivalent to fix a reduced decomposition of $\lam$ in $\Waff/W$, where 
$\Waff$ is the affine Weyl group.
In geometric terms, such a choice is equivalent to fix a Demazure--Hansen--Bott--Samelson 
desingularization $\pi:\hat{\Sigma}(\gamma_\lam) \rightarrow X_\lam$. 
Contou--Carr\`ere \cite{CC} has shown that the points of an affine 
Bott--Samelson variety can be viewed in terms of the affine building $\cjaff$ 
as the set of all galleries 
of the same type as $\gamma_\lam$ starting at the origin ($=G(\co)$). Using this
interpretation in terms of galleries, we show that the birational 
desingularization map identifies the subset of {\it minimal galleries} 
$\hat\cg_\lam\subset \hat{\Sigma}(\gamma_\lam)$ with the
orbit $\cg_\lam\subset X_\lam$.

Using the action of $T$, we get a cellular decomposition of $\hat\Sigma(\g_\lam)$
indexed by the $T$--fixed points. In terms of the interpretation of $\hat\Sigma(\g_\lam)$
as a set of galleries, the $T$--fixed points correspond to the galleries
$\G(\g_\lam)$. This defines a set theoretic map $\hat\Sigma(\g_\lam)\rightarrow \G(\g_\lam)$
sending a point to the center of its cell. This map has in fact also an interpretation in 
terms of buildings.
The affine Tits--building $\cjaff$ is the union of its apartments, and all apartments
are isomorphic to $\ca$. The theory of Tits--buildings provides a retraction 
of the affine building $\cjaff$ onto the affine Coxeter complex $\ca$.
The retraction maps the set $\hat\cg_\lam$ of minimal galleries of type $\gamma_\lam$ in 
$\cjaff$ onto combinatorial galleries in $\ca$ of type $\gamma_\lam$. 
More precisely we get (see  Theorem~\ref{firstpartmaintheorem} and \ref{mvcycledescription}, 
Corollary~\ref{fibrequasiaffine})
\par\vskip 8pt\noindent
{\bf Theorem B.} {\it
The retraction induces
a map $r_{\gamma_\lam}:\hat\cg_\lam\rightarrow \Gamma^+(\gamma_\lam)$ onto the set of all 
positively folded galleries of the same type as $\gamma_\lam$. For such a gallery $\delta$, 
the fibre $r_{\gamma_\lam}^{-1}(\delta)$ is naturally equipped with the structure of an irreducible 
quasi-affine variety, it is the intersection of a Bia{\l}ynicki--Birula cell of $\hat{\Sigma}(\gamma_\lam)$
with $\hat\cg_\lam$. The dimension of the fibre is equal to the combinatorially defined dimension $\dim\de$, 
and $r_{\gamma_\lam}^{-1}(\delta)$ admits a finite decomposition into a union of subvarieties, each being a 
product of $\bc$'s and $\bc^*$'s. In particular, the fibre admits a canonical open and dense subvariety
isomorphic to $\bc^a\times (\bc^*)^b$, where $b = \sharp J^-_{-\infty}(\de)$
(section \ref{cell}) and $a+b=\dim\de$.}
\vskip 8pt
Let $p$ be the bijection $p:\hat\cg_\lam\rightarrow\cg_\lam$, and for
$\delta\in \Gamma^+(\gamma_\lam)$ write $Z(\delta)$ for the closure 
of the fibre in $X_\lam$. Denote by $X_\lam^\mu$ the union of the $Z(\delta)$ for
all $\delta\in  \Gamma^+(\gamma_\lam)$ having $\mu$ as target:
$$
Z(\delta)=\overline{p\big(r_{\gamma_\lam}^{-1}(\delta)}\big)\subset X_\lam,
\quad 
X_\lam^\mu=\bigcup_{\delta\in  \Gamma^+(\gamma_\lam),\ e(\delta)=\mu} Z(\delta).
$$ 
In group theoretic terms, we have $X_\lam^\mu=\overline{U^-(\ck)\mu\cap \cg_\lam}$,
where $U^-\subset B^-$ is the unipotent radical. The special r\^ole played by the 
LS-galleries is that the corresponding $Z(\delta)$ are precisely the
irreducible components of $X_\lam^\mu$ (see Theorem~\ref{mvcycledescription}). 
\par\vskip 8pt\noindent
{\bf Theorem C.} {\it The irreducible components of $X_\lam^\mu$ are given by the $Z(\delta)$
for $\delta$ a LS--gallery, i.e., $X_\lam(\mu)=\bigcup_\delta Z(\delta)$, where $\delta$ runs 
only over all LS--galleries in $\Gamma^+(\gamma_\lam)$ having as target $e(\delta)=\mu$. 
These irreducible components are precisely the MV--cycles.
}
\vskip 8pt
As a consequence, the algebraic cycles in 
$$
MV(\lam)= \{Z(\delta)=
\overline{r_{\gamma_\lam}^{-1}(\delta)}\mid \delta
\hbox{\rm \ a LS--gallery in\ } \Gamma^+(\gamma_\lam)\},
$$
form by \cite{MV} and \cite{Vass} a 
canonical basis of $V(\lam)$, realized as the intersection homology of $X_\lam$. 
So, Theorem C provides the re\-pre\-sentation theoretic 
interpretation of the combinatorial character formula in Theorem~A,
and the construction relates the path model to 
yet another construction of bases for representations  
(see also \cite{Lit3}, \cite{J}, \cite{kashsim}).

It should be very interesting to connect directly the methods presented in this paper with the 
methods (using the moment map) developed by Anderson and Kogan in \cite{An}, \cite{AnK}.

To make the paper as self contained as possible, we recall in section \ref{affineGrassmannian} the definition 
of the {affine Grassmannian} and we explain how we restrict ourselves to the case of a 
{simply connected} semisimple group. The next section deals with the {affine Kac-Moody group} 
associated to $G$ and various incarnations of the affine Grassmannian. The third section concerns 
the {\it building theory}, we recall the definitions of the three buildings one can associate to $G$: 
the {\it spherical building}, the  {\it affine building} and the  {\it building at infinity}, and we 
introduce the {\it retraction centered at $-\infty$}. 

In section \ref{combinatorial gallery}, we define our main objects, the {\it combinatorial galleries}. In 
section \ref{dimensiongallery}, we associate a {dimension} to a positively folded gallery and we 
introduce the {\it LS-galleries}.

To prove the majoration of the dimension of a positively folded gallery, we introduce
in section \ref{rootoperators} the folding operators, which are appropriate analogues
of the root operators for the path model of a representation. We also give a combinatorial 
characterization of the LS-galleries, and  we give the character formula for
$V(\lambda)$ in terms of LS-galleries.

In section \ref{varietiesofgalleries}, we discuss the connection between a
Demazure--Hansen--Bott--Samelson resolution of the Schubert variety $X_\lam$
and galleries of a fixed type in the Tits--building. The last sections are 
devoted to the construction of the stratification of the orbit $\cg_\lam$ 
indexed by the positively folded galleries and the proof of the 
connection between MV-cycles and strata corresponding to LS--galleries.
\vskip 5pt
\noindent
{\bf Acknowledgements:} Both authors thank Henning Haahr Andersen, Edward Frenkel, 
Martin H\"arterich, Shrawan Kumar, Claus Mokler,
Markus Reineke and Guy Rousseau for various helpful discussions. We are especially 
indebted  to Guy Rousseau, who suggested to us the idea to use the retraction $r_{-\infty}$.
The first author would like to thank the CAALT for the support and the 
Aarhus Mathematics Department for the hospitality
during the year 2002/2003.
The second author is happy to thank the MSRI and the Department of Mathematics 
at UC Berkeley for the hospitality and support during the spring semester 2003, where 
part of this article has been worked out.


\section{The Affine Grassmannian}\label{affineGrassmannian} 
Let $G$ be a connected semisimple complex algebraic group. For a commutative $\bc$-algebra ${\mathcal R}$
let $G({\mathcal R})$ be the set of ${\mathcal R}$--rational points of $G$, i.e., the
set of algebra homomorphisms from the coordinate ring $\bc[G]\rightarrow {\mathcal R}$.
Then $G({\mathcal R})$ comes naturally equipped again with a group structure, for example
for $G=SL_n(\bc)$, we can identify $SL_n({\mathcal R})$ with the set of
$n\times n$--matrices with entries in ${\mathcal R}$ and determinant 1. Similarly,
by embedding $G\hookrightarrow SL_n(\bc)$, we can identify $G({\mathcal R})$ naturally
with a subgroup of $SL_n({\mathcal R})$.
 
Denote $\co=\bc[[t]]$ the ring of formal power series in one variable and let $\ck=\bc((t))$ be 
its fraction field, the field of formal Laurent series. Denote $v:\ck^*\rightarrow\bz$
the standard valuation on $\ck$ such that $\co=\{f\in\ck\mid v(f)\ge 0\}$. 
The {\it loop group} $G(\ck)$ is the set of 
$\ck$--valued  points of $G$, we denote by $G(\co)$ its subgroup of $\co$--valued points.
The latter has a decomposition as a semi-direct product $G\semi G^1(\co)$, where
we view $G\subset G(\co)$ as the subgroup of constant loops and $G^1(\co)$ is the
subgroup of elements congruent to the identity modulo $t$. Note that we can describe
$G^1(\co)$ also as the image of $(\Lie G)\otimes_\bc t\bc[[t]]$ via the exponential map.
As a set, the {\it affine grassmannian} $\cg$ is the quotient
$$
\cg=G(\ck)/G(\co).
$$ 
Note that $G(\ck)$ and $\cg$ are {\it ind}--schemes and $G(\co)$ is a group 
scheme (see \cite{BL}, \cite{Ku}, \cite{LS}, \cite{Lu}). There is a model of $\cg$ due to Lusztig
which describes $\cg$  as an increasing union of finite dimensional
complex projective varieties $\cg^{(n)}$, where each $\cg^{(n)}$ is a subvariety of some 
finite dimensional Grassmann variety.

Fix a {\it maximal torus} $T\subset G$ and {\it Borel subgroups} $B,B^-\subset G$ such that 
$B\cap B^-=T$. We denote $\la\cdot,\cdot\ra$ the non--degenerate pairing between 
the {\it character group} $X:=\Mor(T,\bc^*)$ of $T$ and the group $\Xvee:=\Mor(\bc^*,T)$ of 
{\it co-characters}. We identify the lattice $\Xvee$ with the quotient 
$T(\ck)/T(\co)$, so we use the same symbol $\lam$ for the co-character
and the point in $\cg$. 
Let now $p:G'\rightarrow G$ be an isogeny with $G'$ being simply connected. Then
$G'(\co)\simeq G'\semi (G')^1(\co)$, and since $(G')^1(\co)$ and $G^1(\co)$ are
naturally isomorphic, we see that the natural map $p_\co:G'(\co)\rightarrow G(\co)$ is surjective
and has the same kernel as $p$. Let $X'$ and ${X'}^\vee$ be the character group respectively
group of co-characters of $G'$ for a maximal torus $T'\subset G'$ such that $p(T')=T$,
then $p:T'\rightarrow T$ induces an inclusion ${X'}^\vee\hookrightarrow {X}^\vee$. 

The quotient  $\Xvee/{X'}^\vee$ measures the difference between $\cg$ and the affine
grassmannian $\cg'=G'(\ck)/G'(\co)$. In fact, $\cg'$ is connected, and the connected 
components of $\cg$ are indexed by  $\Xvee/{X'}^\vee$. The natural maps 
$p_\ck:G'(\ck)\rightarrow G(\ck)$ and $p_\co:G'(\co)\rightarrow G(\co)$ induce a 
$G'(\ck)$--equivariant inclusion $\cg'\hookrightarrow \cg$, which is an isomorphism onto the 
component of $\cg$ containing the class of $1$. Now $G'(\ck)$ acts via $p_\ck$ on
all of $\cg$, and each connected component is a homogeneous space for $G'(\ck)$,
isomorphic to $G'(\ck)/\cq$ for some parahoric subgroup $\cq$ of $G'(\ck)$ which is
conjugate to $G(\co)$ by an outer automorphism (see \ref{affinebuilding} below, for a 
definition of parahoric subgroups).

\begin{rem}\label{orbits}
\rm
So to study $G(\co)$--orbits on $G(\ck)/G(\co)$ for $G$ semisimple, without loss of generality
we may assume that $G$ is simply connected, but we have to 
investigate more generally $G(\co)$--orbits on $G(\ck)/\cq$ for
all parahoric subgroups $\cq\subset G(\ck)$ conjugate to 
$G(\co)$ by an outer automorphism. 
\end{rem}
\section{The affine Kac-Moody group}\label{KMgroup}
In the following let $G$ be a simply connected semisimple complex algebraic group.  Let again $\co=\bc[[t]]$ denote the ring of formal power 
series in one variable and let $\ck=\bc((t))$ be its fraction field.

The {\it rotation operation} $\gamma:\bc^*\rightarrow \Aut(\ck)$, 
$\gamma(z)\big(f(t)\big)=f(zt)$ gives rise to group automorphisms
$\gamma_G:\bc^*\rightarrow \Aut(G(\ck))$, we denote $\cl(G(\ck))$ 
the semidirect product $\bc^*\semi G(\ck)$.  The rotation operation 
on $\ck$ restricts to an operation of $\co$ and hence we have a natural
subgroup $\cl(G(\co)):=\bc^*\semi G(\co)$ (for this and the following see \cite{Ku},
Chapter 13).

Let $\hat\cl(G)$ be the affine Kac-Moody group associated to
the affine Kac--Moody algebra 
$$
\hat\cl(\Lg)=\Lg\otimes \ck\oplus\bc c\oplus \bc d,
$$
where $0\rightarrow \bc c\rightarrow \Lg\otimes \ck\oplus\bc c\rightarrow
\Lg\otimes \ck\rightarrow 0$ is the universal central extension of the
{\it loop algebra} $ \Lg\otimes \ck$ and $d$ denotes the scaling element.
We have corresponding exact sequences also on the level of groups,
i.e., $\hat\cl(G)$ is a central extension of $\cl(G(\ck))$ 
$$
1\rightarrow\bc^*\rightarrow \hat\cl(G)\mapright{\pi} \cl(G(\ck))\rightarrow 1
$$
(see~\cite{Ku}, Chapter 13).
Denote $\cp_\co\subset \hat\cl(G)$ the ``parabolic'' subgroup $\pi^{-1}(\cl(G(\co)))$.
We have four incarnations of the affine grassmannian:
\begin{equation}\label{fourgrassmann}
\cg=G(\ck)/G(\co)=\cl(G(\ck))/\cl(G(\co))=\hat\cl(G)/\cp_\co = G(\bc[t,t^{-1}])/G(\bc[t]).
\end{equation} 
We consider now the various maximal tori and Weyl groups. Let $N=N_G(T)$
be the normalizer in $G$ of the fixed maximal torus $T\subset G$, we denote
by $W$ the {\it Weyl group} $N/T$ of $G$. Let $N_\ck$ be the subgroup of $G(\ck)$ 
generated by $N$ and $T(\ck)$, 
let $\overline{T}$ be the {\it standard maximal torus}
$\bc^*\times T$ in $\cl(G(\ck))$ and denote by 
$\overline{N}\subset \cl(G(\ck))$ the extension of $N_\ck$.
Finally, let $\ct\subset\hat\cl(G)$ be the standard
maximal torus (such that $\pi(\ct)\subset\overline{T}$), and let $\cn$ 
be its normalizer in $\hat\cl(G)$, then we get three incarnations of the affine Weyl group:
$$
\Waff=N_\ck/T\simeq \overline{N}/\overline{T}\simeq \cn/\ct
$$
Let $ev:G(\co)\rightarrow G$ and  
$ev_\cl:\cl(G(\co))\rightarrow \bc^*\times G$ be the evaluation 
maps at $t=0$, and let $\cb=ev^{-1}(B)$ respectively 
$\cb_\cl=ev_\cl^{-1}(\bc^*\times B)$ be the corresponding Iwahori subgroups,
and let $\hat\cb=\pi^{-1}(B_\cl)\subset \hat\cl(G)$ be the Borel subgroup.
We have the corresponding Bruhat decompositions, and
note that the pair $(\cb,N_\ck)$ is a BN--pair in the loop group $G(\ck)$, 
and $(\hat\cb,\cn)$ is a BN--pair in the affine Kac--Moody group $\hat\cl(G)$.
$$
G(\ck)=\bigcup_{w\in \Waff} \cb w \cb,\quad 
\cl(G(\ck))=\bigcup_{w\in \Waff} \cb_\cl w \cb_\cl,\quad 
\hat\cl(G)=\bigcup_{w\in \Waff} \hat\cb w \hat\cb. 
$$
Using similar identifications as in (\ref{fourgrassmann}), whenever appropriate we may
replace the study of $G(\co)$--orbits on $G(\ck)/\cq$, where
$\cq$ is a parahoric subgroup of $G(\ck)$ conjugate to 
$G(\co)$ by an outer automorphism of $G(\ck)$, by 
the study of $\cp_\co$--orbits on $\hat\cl(G)/\hat\cq$, where 
$\hat\cq\subset \hat\cl(G)$ is a parabolic subgroup of $\hat\cl(G)$,
conjugate to $\cp_\co$ by a diagram automorphism of $\hat\cl(G)$. 
 These orbits correspond exactly to each other
because the kernel of $\pi$ acts trivially on  $\hat\cl(G)/\hat\cq=\cl(G(\ck))/\pi(\hat\cq)$, 
so $\cl(G(\co))$ acts naturally on $\hat\cl(G)/\cq$, and the $\cl(G(\co))$--orbits 
and the $G(\co)$--orbits in $\hat\cl(G)/\hat\cq$ obviously coincide since 
by the Bruhat decomposition the orbits are
parameterized by the same ${T}$-- respectively
$\overline{T}$-- respectively  $\ct$--fixed points in 
$G(\ck)/\cq=\hat\cl(G)/\hat\cq$.


\section{Buildings, roots and characters}\label{buildingsrootsandcharacters}

From an abstract point of view, a building $\cj$ is a simplicial complex covered by some 
subcomplexes, called apartments, such that, each apartment is a Coxeter complex, any 
two simplices of $\cj$ are always contained in an apartment, and given two apartments 
$A$ and $A'$ with a common maximal simplex, then there is an isomorphism 
$A\to A'$ fixing every simplex in $A \cap A'$ (see for example \cite{B} or \cite{T74}).

These properties can be used to define retractions $r:\cj\rightarrow A$. We will use the 
building theory later because one can identify a $G(\co)$--orbit in 
$G(\ck)/G(\co)$ with a subset of the affine building (see \ref{affinebuilding}) associated 
to $G$ (not needed to be simply connected). In particular, the retraction of the building onto an apartment will be the
main tool to identify certain strata of a $G(\co)$--orbit
with combinatorial objects in the apartment. Note that the apartment 
is in our case nothing else then the Coxeter complex of the affine Weyl group.

In this section we fix some notation and recall for the convenience
of the reader some general facts on building theory. As 
references we suggest  \cite{B}, \cite{BT}, \cite{R}  and/or \cite{T74}.

\subsection{Roots and characters}
We denote $\la\cdot,\cdot\ra$ the non--degenerate pairing between 
the character group $X:=\Mor(T,\bc^*)$ of $T$ and its group $\Xvee:=\Mor(\bc^*,T)$ of
co-characters. Let $\Phi\subset X$ be the root system of the pair 
$(G,T)$, and, corresponding to the choice of $B$, denote $\Phi^+$ the set of positive roots, 
let $\Delta=\{\al_1,\ldots,\al_n\}$ be the set of simple roots, and let $\rho$ be half the 
sum of the positive roots.
Let $\Phivee\subset \Xvee$ be the dual root system, together with a bijection 
$\Delta\rightarrow \Delvee$, $\al\mapsto\alvee$.  We denote $\Rveep$ the 
submonoide of the coroot lattice $\Rvee$ generated by the positive coroots $\Phiveep$.
We define on $\Xvee$ a partial order by setting $\lam\succ \nu\Leftrightarrow \lam-\nu\in\Rveep$, 
and let $\Xveep$ be the cone of 
dominant co-characters.
$$\Xveep:=\{\lam\in\Xvee \mid \la\lam,\al\ra\ge 0\,\forall\,\al\in\Phi^+ \}$$

In the following, we will deal with three buildings associated to the simply connected 
complex semisimple group $G$. The first building we consider is the
\subsection{Spherical Building}\label{sphericalbuilding}
Let us denote by $\cj^s$ the set of all the parabolic subgroups of $G$. 
The opposite relation of the inclusion between parabolic subgroups 
endows this set with a structure of simplicial complex. The maximal 
simplices given by the Borel subgroups will be called ``spherical" 
chambers, the others, simplices or ``spherical" faces. We will often 
drop the term ``spherical" when no confusion may arise. 

For every maximal torus $T$ in $G$, we define an apartment in $\cj^s$, 
by letting $\A^s(T)$ be the set of all the parabolic subgroups that contain $T$. 
Together with those apartments, the simplicial complex $\cj^s$ is a building 
(see for example \cite{T74} or \cite{B}). The apartments are isomorphic to the 
Coxeter complex $C(W,S)$ defined by the Coxeter presentation of the Weyl group $W$. 
The existence of the retractions in the spherical building is equivalent to the Bruhat 
decompositions of the group $G$. 
\begin{exam}\label{exampleBruhat}
\rm
Denote by $\Lc_f$ the chamber corresponding to the Borel subgroup $B$.
Recall that, as a set, we can identify the set of chambers in the building
with the flag variety $G/B$. The retraction $r_{\Lc_f}: \cj^s\to\A^s$ 
onto the
apartment $\A^s$ of center $\Lc_f$ is defined as follows:
Given a chamber $F_{B'}$ in the building associated to a Borel subgroup 
$B'$, by the
Bruhat decomposition we can find $b\in B$ and $w\in W$ such that $B' = 
bwB/B$ in
$G/B$, and the image by the retraction is the $T$-fixed point $wB/B$ in 
the Schubert
cell $BwB/B$. In the language of buildings we set $r_{\Lc_f}(F_{B'}):= 
F_{wBw^{-1}}$.
This map is actually defined on any face since there exists only one 
parabolic subgroup of a
given type containing a given Borel subgroup.
\end{exam}

\subsection{Affine Building}\label{affinebuilding}
We can identify in a similar way the affine building $\cj^{\mathfrak a}$,
as a set, with the set of all the parahoric subgroups of $G(\ck)$, that is the set of 
all the subgroups of $G(\ck)$ that contain a conjugate of $\cb$. The construction 
below differs from the one above, but, with the appropriate changes, it
would apply in the same way also to the spherical building.

Set $\tta:=X^\vee\otimes_\bz\br$, the affine Weyl group
$\Waff$ acts on $\tta$ as an {\it affine reflection group}. Let $\cha\subset \tta$
be the set of affine reflection hyperplanes for the action of $\Waff$ (i.e., the affine reflections 
with respect to these hyperplanes generate $\Waff$ and the set $\cha$ is stable under
the action of $\Waff$). The connected components of $\tta-\cha$ are called
{\it open alcoves}, the closure of such a component is called a 
{\it closed alcove} or just {\it alcove}.

The hyperplanes in $\cha$ are all of the form 
$\ch_{\beta,m}=\{a\in\tta\mid \langle a,\beta\rangle=m\}$
for some positive root $\beta\in\Phi^+$ and $m\in\bz$, we denote $s_{\beta,m}$ the 
corresponding affine reflection. The associated closed affine halfspaces are
denoted $\ch_{\beta,m}^+=\{a\in\tta\mid \langle a,\beta\rangle\ge m\}$, 
resp. $\ch_{\beta,m}^-=\{a\in\tta\mid \langle a,\beta\rangle\le m\}$. 
We use the notation $\ch_{\beta,m}^{+,o}$,  
respectively $\ch_{\beta,m}^{-,o}$, for the corresponding open affine halfspaces.

\begin{dfn}\rm
By a {\it face} $F$ we mean a subset of $\tta$ obtained as the intersection 
of closed affine halfspaces and affine hyperplanes, the intersection running 
over all pairs $(\beta,m)$, $\beta\in\Phi^+$, $m\in\bz$. By the corresponding
 {\it open face} $F^o$ we mean the subset of $F$ obtained when replacing
 the closed affine halfspaces in the definition of $F$ by the corresponding
 open affine halfspaces. 
 \end{dfn}
The open faces define a partition of $\tta$: we call two elements
$x,y\in \tta$ equivalent if for all pairs $(\beta,m)$, $x$ and $y$ are both 
in either $\ch_{\beta,m}$, or $\ch^{+,o}_{\beta,m}$, or $\ch^{-,o}_{\beta,m}$. 
The equivalence classes of this relation are the open faces $F^o$, and the closed 
face $F$ is the closure of $F^o$ in the affine subspace $\langle F^o\rangle_{\rm aff}$ 
spanned by $F^o$. 

We call  $\langle F^o\rangle_{\rm aff}=\langle F\rangle_{\rm aff}$ the 
{\it support} of the (open) face, the {\it dimension} of the face is the dimension of its support. 
So the alcoves are the faces of maximal dimension. Let $F$, $C$ be two faces. We say that $F$ is a face of 
$C$ if $F$ is defined by changing some inequalities in the definition of $C$ into equalities. 
A {\it wall of an alcove} is the support of a codimension one face. In general, instead of
the term hyperplane we use often the term {\it wall}, which is more common in the language of buildings.
The following is well-known:

\begin{thm}
The affine Weyl group $\Waff$ acts simply transitive on the set of all alcoves.
The {\it fundamental alcove} $\Delta_f=\{\nu\in \tta\mid  
0\le\langle x,\be\rangle\le 1, \forall\,\be\in\Phi^+\}$ is a fundamental
domain for the action, and $\Waff$ is generated by the affine reflections 
$$
\Saff=\{s_{\be,m}\mid \ch_{\beta,m} \hbox{\rm\ is a wall for\ }\Delta_f\}
$$
\end{thm}
Denote by $S\subset \Saff$ the set of reflections 
$S=\{s_\beta:=s_{\be,0}\mid \ch_{\beta,0} \hbox{\rm\ is a wall for\ }\Delta_f\}$; we can describe
$S$ also as $\Saff(0)=\{s_{\beta,m}\in\Saff\mid  0\subset \ch_{\beta,m}\}$.
More generally, for a face $F$ of $\De_f$ let 
$$
\Saff(F)=\{s_{\beta,m}\in\Saff\mid  F\subset \ch_{\beta,m}\}.
$$
We call $\Saff(F)$ the {\it type of $F$}, so $\Saff(0)=S$ and $\Saff(\De_f)=\emptyset$.
For an {\it arbitrary face $F$ its type} is defined as the type of $F'$, where
$F'$ is the unique face of $\Delta_f$ such that $F=w(F')$ for some $w\in\Waff$.

\begin{rem}\rm
The alcoves are actually the chambers of the {\it Coxeter complex} $\tta^{\mathfrak{a}}$ 
associated to the Coxeter group $(\Waff,\Saff)$, but since we look at the same time
also at the spherical complex, to not confuse the ``affine'' and the spherical chambers,
we prefer the term alcove for the faces of maximal dimension. 
\end{rem}
We view the (spherical) Weyl group $W\subset\Waff$ as the subgroup generated by the 
reflections in $S$. By an {\it open chamber} we mean always an open spherical chamber, 
i.e., a connected component of $\tta-\bigcup_{\beta\in\Phi^+} \ch_{\beta,0}$, and 
a {\it chamber} is the closure of such a connected component. We have the
{\it dominant chamber $\Lc_f$} corresponding to the choice of the Borel subgroup 
$B$ and the {\it anti-dominant chamber} $-\Lc_f$.
\begin{dfn}\rm
A {\it sector} $\Ls$ in $\tta$  is a $\Waff$--translate of a chamber. Two sectors $\Ls,\Ls'$
{\it are called equivalent} if there exists a third sector $\Ls''$ in the intersection: 
$\Ls''\subset\Ls\cap\Ls'$.
\end{dfn}
\begin{lem}
The equivalence classes of sectors are in one-to-one correspondence with the spherical
chambers, i.e., in every class there exists a unique spherical chamber.
\end{lem}
For a root $\beta\in\Phi$ let $\Lg_\beta\subset\Lg$ be the root subspace of the complex
Lie algebra $\Lg=\Lie G$ and fix a generator $X_\beta$. Via the exponential map 
$\exp:\bc \rightarrow G$, $x\mapsto \exp(xX_\beta)$ we get a one-dimensional
unipotent subgroup $U_\beta$ normalized by $T$.
\begin{dfn}\rm
For a real number $r$ let $U_{\beta,r}\subset U_\beta(\ck)$ be the unipotent subgroup
$$
U_{\beta,r}=\{1\}\cup \Big\{exp(X_\beta\otimes f )\mid f\in\ck^*, v(f)\ge r \Big\}.
$$
For a non-empty subset $\Omega\subset \tta$ let 
$\ell_\beta(\Omega)=-\inf_{x\in \Omega} \beta(x)$. We attach to 
$\Omega$ a subgroup of $G(\ck)$ by setting
\begin{equation}\label{uomega}
U_\Omega:=\langle U_{\beta,\ell_\beta(\Omega)}\mid\beta\in\Phi\rangle.
\end{equation}
\end{dfn}\rm
\begin{exam}\rm
If $\Omega=0$, then we have $U_0=G(\co)$, and if $w\in \Waff$, then 
$U_{w(\Omega)}=n_wU_\Omega n_w^{-1}$ for any representative $n_w\in N(\ck)$ of $w$. 
\end{exam}
To define the affine building $\cj^{\mathfrak a}$,
let $\sim$ be the relation on $G(\ck)\times\tta$ defined by:
$$
(g,x)\sim (h,y)\quad\hbox{\rm if }\exists\,n\in {N(\ck)}\ \hbox{\rm such that\ }
nx=y\ \hbox{\rm and\ }g^{-1}hn\in U_x.
$$
\begin{dfn}\rm
The {\it affine building} $\cj^{\mathfrak a}:=G(\ck)\times\tta/\sim$ associated to $G$ is the
quotient of $G(\ck)\times\tta$ by ``$\sim$''.
The building $\cj^{\mathfrak a}$ comes naturally equipped with a $G(\ck)$--action
$g\cdot(h,y):=(gh,y)$ for $g\in G(\ck)$ and $(h,y)\in \cj^{\mathfrak a}$.
\end{dfn}
The map $\tta\rightarrow \cj^{\mathfrak a}$, $x\mapsto(1,x)$
is injective and $N(\ck)$ equivariant, we will identify in the following 
$\tta$ with its image in $\cj^{\mathfrak a}$.

The stabilizer $P_{(g,y)}$ of any point is a parahoric subgroup. In fact, 
for $x\in \tta$ it is the parahoric subgroup generated by the stabilizer 
$N_x=\{n\in N(\ck)\mid nx=x\}$ and $U_x$:
$$
P_{x}=\langle U_x,\ N_x \rangle.
$$
\begin{exam}\rm
If $x\in \Delta_f^o$ is in the open face of the fundamental alcove $\Delta_f$, then 
$P_x=\cb$. More generally, let $F$ be a face of the fundamental alcove and denote
by $\Waff(F)\subset \Waff$ the subgroup generated by its type $\Saff(F)$.
If $x\in F^o$ is an element in the corresponding open face,
then $P_x=P_F$ is the corresponding {\it standard parahoric
subgroup of type $F$}:
$$
P_F:=\bigcup_{w\in \Waff(F)} \cb w\cb.
$$
\end{exam}
\begin{dfn}\rm
The subsets of $\cj^{\mathfrak a}$ of the form $g\tta$ are called {\it apartments}.
\end{dfn}
\begin{lem}
For any $g\in G(\ck)$ and $x\in \tta\cap g^{-1}\tta$, there exists a $n\in N(\ck)$ such that 
$gx=nx$, or, in other words, $G(\ck) x\cap \tta= N(\ck)x$ for $x\in\tta$.
\end{lem}
So we can extend the partition of $\tta$ into open faces
into a partition of $\cj^{\mathfrak a}$ by calling a subset $F^o$
of $\cj^{\mathfrak a}$ an open face if it is of the form $g{F'}^o$ for some
face $F'\in\tta$. Similarly, we call a subset of $\cj^{\mathfrak a}$ a face, an alcove,
a sector etc. if it is of the form $gF$, $g\Delta$, $g\Ls$ etc. for some
face, alcove, sector etc. in $\tta$. Note that for any non-empty set 
$\Omega\subset \tta$ the group $U_\Omega$ acts transitively on the set of 
all apartments which contain $\Omega$.

We have a Bruhat decomposition $G(\ck)=U_xN(\ck)U_y$
for any pair of elements $x,y\in\tta$. Another important fact
from the theory of buildings is that any two faces
or any two alcoves are contained in a common apartment, and if 
$\Ls,\Ls'$ are two sectors, then there exist subsectors $\Ls_1\subset \Ls$
and $\Ls_1'\subset \Ls'$ such that $\Ls_1$ and $\Ls_1'$ are contained in
a common apartment, and if $\Ls$ is a sector and $\Delta$ is an alcove, then there
exists a subsector $\Ls_1\subset \Ls$ such that $\Delta$ and $\Ls_1$ are in a common apartment.

We see that on the level of sets the affine building $\cj^{\mathfrak a}$ is 
the disjoint union of $G(\ck)/P_F\times F^o$'s, where the $F$'s are running through 
the set of faces of the fundamental alcove.

\begin{rem}\rm
The action of $G(K)$ extends to an action of the Kac-Moody group $\hat\cl(G)$ on $\cj^{\mathfrak a}$ such that the stabilizer of $\De_f$ is the Borel subgroup $\hat\cb$ and the center acts trivially.
\end{rem}

\subsection{The Building at Infinity and the Retraction $\rinf$}\label{inftyandretraction}

The last building we introduce is called the {\it spherical building at infinity}, we 
will denote it by $\cj^{\infty}$. We refer to \cite{B} or \cite{R} for a precise definition, 
we recall quickly its construction and the properties which we will need. 

The apartments for $\cj^{\infty}$ are the same as for $\cj^{\mathfrak{a}}$, only the structure 
of the complex is now different. The chambers of the complex in an apartment
are the equivalence classes of sectors. So the structure is similar to that of the
spherical complex $\cj^s$, only that one does not have anymore a preferred
vertex for the ``Weyl chambers''. Still, the apartments of $\cj^{\infty}$ are 
isomorphic to the Coxeter complex $C(W,S)$ of the Weyl group $W$.
Since the equivalence classes of sectors are in one-to-one correspondence 
with the spherical chambers, it makes sense to use the following notation:
\begin{dfn}
The equivalence class of sectors of the anti-dominant Weyl chamber $-\Lc_f$
is called the anti-dominant chamber at $-\infty$ and is denoted $\Lc_{-\infty}$.
\end{dfn}
For any alcove $\De$ in the apartment $\tta \hookrightarrow \cj^{\mathfrak a}$, one can define a 
chamber complex map $r_{\De,\tta} : \cj^{\mathfrak a}\to \tta$, called the retraction onto 
$\tta$ of center $\De$. The retraction has the following two properties, which in addition
characterize the map uniquely: the first property is that for any face $F$ of 
$\De$, $r_{\De,\tta}^{-1}(F) = \set{F}$. 

To recall the second property, we define a distance on the set of alcoves. 
Two alcoves are called {\it adjacent} if they have a common codimension one face.
A {\it gallery of alcoves of length $r$} is a sequence $(\Delta_0,\ldots,\Delta_r)$ of alcoves 
such that $\Delta_i$ and $\Delta_{i-1}$ are adjacent for $1\le i\le r$.                        
The distance $d(\Delta,\Delta')$ of two alcoves in $\cj^{\mathfrak a}$ is the minimal
length of a gallery $(\Delta_0,\ldots,\Delta_r)$ joining $\Delta$ and $\Delta'$ (i.e.,
$\Delta_0=\Delta$ and $\Delta_r=\Delta'$). 

A gallery of alcoves $(\Delta_0,\ldots,\Delta_r)$ is called {\it minimal} if the length
is equal to the distance $d(\Delta_0,\Delta_r)$ between the first and the last alcove.
Recall (\cite{BT}, section 2.3.6), two alcoves lie always in an apartment containing both, and any apartment
that contains both also contains all minimal galleries joining the two.

The second important property of the retraction is that for any alcove $\De'$, 
$d(\De, \De')= d\big(\De, r_{\De,\tta}(\De')\big )$, so the map preserves the distance
from its center. Further, $r_{\De,\tta}$ restricts to an isomorphism of chamber complexes 
$g\tta\simeq\tta$, for any apartment $g\tta$, $g\in G(\ck)$.

The image $ r_{\De,\tta}(\De')$ of an alcove $\De'$ depends of course on the alcove
$\De$, but the image of $\Delta'$ becomes in fact stable if $\De$ is only ``far away'' enough
from $\Delta'$. To be more precise, consider an alcove $\De'\in \cj^{\mathfrak a}$ and let
$-\Lc_f$ be the anti-dominant chamber. Then there exists a sector 
$\Ls$ equivalent to $-\Lc_f$
such that $\Ls$ and $\Delta'$ lie in a common apartment $g\tta$.  Let $\Delta$
be an alcove in the sector and recall that $r_{\De,\tta} $ restricts to an isomorphism 
of chamber complexes  $g\tta\simeq\tta$, which fixes of course the common sector $\Ls$.
It follows (see for example, \cite{R} $\S$9.4) that $r_{\De,\tta}(\De')$ is independent
of the choice of $\Delta\subset\Ls$. 
\begin{dfn}\rm
The map $r_{-\infty}:\cj^{\mathfrak a}\rightarrow\tta$, defined by 
$\rinf(\De')=r_{\De,\tta}(\De')$ for some alcove $\De\subset\Ls$, where
$\Ls$ is a sector, 
equivalent to $-\Lc_f$, contained in a common apartment with 
$\De'$, is called the {\it retraction of center} $-\infty$.
\end{dfn}
The retraction $r_{\De,\tta}$ can also be expressed in terms of
the action of the group $U_\De$ (see~(\ref{uomega})). Recall that the latter
operates transitively on the set of apartments containing $\Delta$, and the
retraction is in fact the projection $r_\Delta:u\tta\rightarrow \tta$
that maps a face $uF$, $u\in U_\Delta$, $F\subset \tta$, onto $F$.
So the fibres of the retraction $r_\Delta$ are exactly the $U_\De$--orbits.

Consider the retraction $\rinf$ of center $-\infty$. For a face
$F\in \tta$ and a face $F'$ in the fibre, we can find a subsector
$\Ls\subset -\Lc_f$ such that for all $\De\in\Ls$ there exists an $u\in U_\Delta$ 
so that $F'=uF$. This means for all $N\gg 0$ we can choose an alcove
in the anti-dominant chamber such that $\De\subset \Ls$ and 
$\ell_\beta(\De)>N$ for all positive roots. This means (compare \cite{BT}, $\S$6 and $\S$7)
that we can actually choose $u\in U^-(\ck)$, where $U^-\subset B^-$ is the unipotent
radical of the opposite Borel subgroup $B^-$. Summarizing we have the 
following description of the fibres of $\rinf$ in terms of the group
operation on $\cj^{\mathfrak a}$:

\begin{prop}\label{fibresofrinf}
The fibres of $\rinf:\cj^{\mathfrak a}\rightarrow\tta$ are the 
$U^-(\ck)$--orbits on $\cj^{\mathfrak a}$.
\end{prop}


\section{Generalized Galleries}\label{combinatorial gallery}
In the next three sections, we take $G$ to be a semisimple group of any type. We first need a more general version of a gallery as the one in section~\ref{buildingsrootsandcharacters}, 
we will essentially follow \cite{CC}.
\begin{dfn}\rm
A {\it generalized gallery in the affine building} is a sequence of faces $\g$ in $\cj^{\mathfrak a}$
$$
\g = (\G'_0\subset \Gamma_0\supset \Gamma_1'\subset\cdots\subset\Gamma_{j-1}\supset
\Gamma_j'\subset  \Gamma_{j}\supset\cdots \supset\Gamma_p' \subset \Gamma_{p}\supset \G'_{p+1}),
$$ such that
\begin{itemize}
\item the first and the last faces, $\G'_0$ and $\G'_{p+1}$, in other words, the {\it source} and the {\it target} of $\g$, 
are vertices of $\cj^{\mathfrak a}$,
\item the $\Gamma_j$'s are faces, all of the same dimension,
\item the $\Gamma_j'$'s, for $j=1,\ldots p$, are faces of two consecutive faces, of relative codimension one.
\end{itemize}
If such a gallery is contained in the apartment $\tta$, it will be called a {\it combinatorial gallery}. 
\end{dfn}

\begin{rem}
\rm
The last two conditions are not necessary to define the gallery model, but we will keep them throughout in the following.
\end{rem}

For any subset $\Omega$ and any face $F$ contained in an apartment $A$ of the affine 
building $\mathcal J^{\mathfrak a}$, we say that a wall $\ch$  
{\it separates $\Omega$ and $F$} if  $\Omega$
is contained in the corresponding closed half space and $F$ is a subset of the opposite
open half space. Let $E$ and $F$ be two faces in the building. In any apartment $A$ 
containing both of them, there exists a finite number of walls that separate $E$ and $F$. 
We denote this set by $\M(E, F)$. Note, if $E$ is not an alcove, all the walls containing 
$E$ but not $F$ belong to $\M(E,F)$. 

Let $\De\supset E$ be an alcove in $A$, then, by \cite{T74}, Proposition 2.29, 
there exists a unique alcove in $A$, denoted $proj_F(\De)$, such that any face of the 
convex hull of $\De$ and $F$ containing $F$ is contained in $proj_F(\De)$.

\begin{dfn}\rm
The alcove $\De$ is said to be {\it at maximal distance to $F$} if the 
length of the minimal alcoves galleries between $\De$ and $proj_F(\De)$ is $\sharp \M(E,F)$. 
\end{dfn} 

Note, any of these minimal galleries will cross only once each wall of $\M(E,F)$.
Such an alcove $\De\supset E$ at maximal distance to $F$ is not uniquely determined.
 
 \begin{lem}\label{Stab}
The alcoves $\De\supset E$ in $A$ at maximal distance to $F$ are conjugate under the
stabilizer of  $E\cup F$.
\end{lem}

\proof
Let $\De'$ be another alcove in $A$, the fixed apartment which contains the two faces $E$ and $F$.
Suppose $\De'\supset E$ is at maximal distance to $F$ and $\Delta$ and $\Delta'$ are adjacent, then
the support $\ch$ of the common face contains $E$. Now for any wall $\ch'\not=\ch$, 
the relative position of $\Delta$ and $\Delta'$ is the same, i.e., they lie in the
same closed halfspace. Since  they lie on different sides of $\ch$, the two can be at maximal 
distance to $F$ at the same time only if $F\subset \ch$. In the general case, choose a minimal gallery
$(\Delta_0=\Delta,\ldots,\Delta_r=\Delta')$. Since $E\subset\Delta,\Delta'$, one concludes
by the minimality of the gallery that $ E\subset \Delta_j$ for all $j$. Further,
for all walls in $\M(E,F)$, $\Delta$ and $\Delta'$ lie in the same closed halfspace,
and hence, by the minimality of the gallery, so do all $\Delta_j$. But this implies
that all the $\Delta_j$ are at maximal distance to $F$, which finishes the proof
since the $\Delta_j$ are pairwise adjacent
\endpf

We have also to generalize the notion of a minimal gallery. Roughly speaking, a generalized
gallery is called minimal if it can be embedded in a minimal gallery of alcoves. More precisely:
\begin{dfn}\label{gengallery}\rm
A generalized gallery 
$$
\g = (F_f = \G'_0\subset \Gamma_0\supset \Gamma_1'\subset\cdots\subset\Gamma_{j-1}
\supset\Gamma_j'\subset  \Gamma_{j}\supset\cdots \supset\Gamma_p' \subset \Gamma_{p}\supset F)
$$ in the affine building $\cj^{\mathfrak a}$ is called {\it minimal} 
if the following holds: let $\De\supset F_f$ be an alcove at maximal distance to $F$, then

{\it a)\/} there exists a minimal gallery of alcoves between $\De$ and $proj_F(\De)$
$$
\mu = (\De=\De_1^0,...,\De_{q_0}^0,...,\De_1^j,...,\De_{q_j}^j,....,\De_1^p,...,\De_{q_p}^p= proj_F(\De)),
$$ 
and for all $j$ the alcoves $\De_1^j,...,\De_{q_j}^j$ contain $\Gamma'_j$. 

{\it b)\/} let $A$ be an apartment containing $\g$, and denote $\ch^A$ the set of affine hyperplanes
in $A$. The set $\M(F_f,F)$ is the disjoint union of the sets 
$$
\sh_j:=\{\ch\in \ch^A\mid \G'_j\subset\ch,\ \G_{j}\not\subset \ch\},\ j=0,\ldots,p.
$$
\end{dfn}
\begin{rem}\rm
A minimal gallery $\g$ is contained in an apartment by condition {\it a)\/}, and an apartment
containing $\De$ and $proj_F(\De))$ contains also all minimal galleries between the two alcoves.
One checks easily that if {\it b)\/} holds for one apartment, then it holds for all apartments $A\supset\g$.
\end{rem}
\begin{rem}\label{minismallgall}
\rm
We will also need to consider minimal galleries of the form 
$(E\supset F'\subset F)$ in an appartment $A$. Such a gallery is 
called {\it minimal} if the set of walls $\M(E,F)$ is exactly the set of 
walls that contain $F'$ but not $F$.

If $\gamma$ is a minimal gallery (in the  sense of Definition~\ref{gengallery}), then 
all $(\Gamma_{j-1}\supset \Gamma'_{j}\subset\Gamma_j)$, $j=1,\ldots,p$ are minimal.
\end{rem}

\begin{rem}\label{equivminisgall}\rm
Using the alcoves at maximal distance, we can reformulate the condition of 
$(E\supset F'\subset F)$ to be minimal
as follows: Let $\mathcal H(F',F)$ be the set of walls that contain $F'$ but not $F$. Then
$(E\supset F'\subset F)$ is a minimal gallery in $A$ if for any alcove $\De\supset E$ in $A$
at maximal distance to $F$, any minimal alcoves gallery between $\De$ and $proj_F(\De)$ has length
$\sharp \mathcal H(F',F)$.
\end{rem}
Our main objects will be minimal galleries having the vertex corresponding to $G(\co)$ as source. 
Since all apartments containing this vertex are conjugate under the action of $G(\co)$, the following 
definitions make sense for arbitrary galleries contained in an apartment (recall, if $\g$ is minimal
then it is contained in at least one apartment). Our aim is to express the condition on the minimality 
in terms of the walls crossed by the gallery.

Let $\lam\in \Xvee$ be a co-character and denote $\ch_\lam:=\bigcap_{\la \lam,\al\ra = 0} \ch_{\al,0}$
the intersection of all the hyperplanes corresponding to the positive roots orthogonal  to $\lam$. 
Let $F_{f}$ be the face (of type $S$) of $\tta$ corresponding to the origin of $\tta$ and let $F_{\lam}$ 
be the face corresponding to $\lam$. 
\begin{dfn}\rm
A {\it combinatorial gallery joining $0$ with $\lam$} is a generalized gallery $\gamma$ in $\tta$ that starts at 
$F_f$ and ends in $F_\lam$:
$$
\g = (F_f = \G'_0\subset \Gamma_0\supset \Gamma_1'\subset\cdots\subset\Gamma_{j-1}\supset
\Gamma_j'\subset  \Gamma_{j}\supset\cdots \supset\Gamma_p' \subset \Gamma_{p}\supset F_{\lam}).
$$
such that the dimension of the faces $\G_j$'s is always equal to $\dim \ch_\lam$.
The latter will be called the {\it large faces} of the gallery $\g$ and the $\G'_j$'s 
the {\it small faces} of $\g$.  
\end{dfn}

\begin{lem}
A combinatorial gallery $\g=(F_f = \G'_0\subset \Gamma_0\supset\cdots \supset\Gamma_p' \subset 
\Gamma_{p}\supset F_{\lam})$  joining $0$ with $\lam$ is minimal if and only if $\g\subset\ch_\lam$ and $p$ is 
minimal in the following sense:\par
$\bullet$ for $j = 0,1,..., p$, let $\sh_j$ be the set of all the affine hyperplanes $\ch\in \cha$ 
such that $\G'_j\subset\ch$ and  $\G_{j}\not\subset \ch$, then the sets $\sh_j$ are pairwise 
distinct and $\cup_{j\in\set{0,...,p-1}}\sh_j = \M(F_f,F_\lam)$, the 
set of all the affine hyperplanes separating any alcove $\De\supset F_f$ at 
maximal distance to $F_\lam$ and the latter.
\end{lem}
\proof
Suppose $\g$ is minimal, then $0,\lam\in\ch_\lam$ implies $\g\subset\ch_\lam$ by minimality, and
the second part of the condition follows from the definition of a minimal gallery.

For the reverse implication, we construct the minimal gallery of alcoves
inductively. We start with an alcove $\Delta=\Delta_1^0\subset F_f$ at maximal distance
to $\Gamma_0$, set $\Delta_{q_0}^0=proj_{\Gamma_0}(\Delta)$. If $p=0$, then the corresponding minimal gallery
joining $\Delta_1^0$ and $\Delta_{q_0}^0$ has length $\sharp\sh_0$ and is already the desired
minimal gallery of alcoves. If $p\ge 1$, note that $\Delta_{q_0}^0$ is at maximal distance to 
$\Gamma_1$: a wall separating $\Delta_1^0$ and $\Gamma_1$ is either
an element of $\sh_0$ or $\sh_1$, so by construction the walls in $\sh_1$ 
separate $\Delta_{q_0}^0$ and $\Gamma_1$. But this implies that the length of
a minimal gallery joining $\Delta_{q_0}^0$ and $\Delta_{q_1}^1=proj_{\G_1}(\Delta_{q_0}^0)$
is $\sharp\sh_1$. By repeating the procedure,  we obtain the desired minimal 
gallery containing $\g$.  
\endpf

Let $\lam\in \Xveep$ be a dominant co-character and
let $\g_{\lam}$ denote a minimal combinatorial gallery joining $F_f$ with $F_{\lam}$ 
$$
\g_{\lam} = (F_f\subset \Upsilon_0\supset \Upsilon_1'\subset\cdots\subset\Upsilon_{j-1}\supset
\Upsilon_j'\subset  \Upsilon_{j}\supset\cdots \supset\Upsilon_p' \subset \Upsilon_{p}\supset F_{\lam}).
$$ 
Because of the minimality assumption, all the faces of $\g_\lam$ are contained in the 
(spherical) dominant chamber $\Lc_f$, and $\Upsilon_0$ is a face of the fundamental alcove $\De_f$. 

\begin{dfn}\rm
The {\it gallery of types} associated to $\g_{\lam}$ is the list of the types of the faces in the
gallery above: 
$$ 
t_{\g_{\lam}} = type(\g_{\lam}) = (S = t'_0\supset t_0\subset t'_1\supset\cdots\supset
t_{j-1}\subset t'_j \supset t_j \subset\cdots \subset t'_p\supset t_p\subset t_{\lam}),
$$ 
where $t_{\lam}$ is the type of the face $F_{\lam}$ and $t_j$ is the type of the face $\Upsilon_j$, 
$t'_j$ the type of $\Upsilon'_j$.
\end{dfn}

\begin{rem}\label{minismallgalllocal}\rm
Fix a minimal combinatorial gallery $\g_\lam$ 
joining $F_f$ with $F_{\lam}$, and consider the set of all generalized  galleries
$\de$ in the affine building $\cj^{\mathfrak a}$ of type $t_{\g_\lam}$. 
Let $\g = (F_f = \G'_0\subset \Gamma_0\supset\cdots \supset\Gamma_p' 
\subset \Gamma_{p}\supset F_{\nu})$ be such a gallery, and assume
that all $(\Gamma_{j-1}\supset \Gamma'_{j}\subset\Gamma_j)$, $j=1,\ldots,p$ are minimal.
Then one sees easily that $\gamma$ is a minimal gallery, the fixed type comes
from a minimal gallery and hence forces the gallery to be minimal.
\end{rem}


Let $\Gamma (\g_\lam)$ be the set of all the combinatorial galleries of 
type $t_{\g_{\lam}}$ and of source $F_f$, i.e., an element 
$\delta\in \Gamma (\g_\lam)$ is a gallery starting at $F_f$:
$$
\delta=(F_f\subset \Sigma_0\supset {\Sigma_1}'\subset\cdots 
\supset{\Sigma}'_p \subset {\Sigma}_{p}\supset F_{\nu}),
$$
such that $\Sigma_j$ is a face of type $t_j$ of $\tta$ and ${\Sigma_j}'$ is of type $t'_j$ and 
$F_\nu$ is of the same type as $F_\lam$. In addition, such a sequence of types gives rise to a 
sequence of subgroups of $\Waff$ 
$$
W'_0\supset W_0 \subset W'_1\supset  \cdots \supset W_{j-1}\subset W'_j 
\supset W_j\subset  \cdots \subset W'_p \supset W_p, 
$$ where for $j = 0,...,p$, $W'_j$ respectively $W_j$ is the Coxeter subgroup of $\Waff$ generated by 
the reflections in $t'_j\subset\Saff$ respectively in $t_j\subset\Saff$. 

These groups are also the stabilizers in $\Waff$ of the faces of the fundamental alcove of the 
corresponding type. This provides us with another way to ``encode'' a gallery of type $\lam$:
The stabilizer of the $F_f$ in $\Waff$ is $W = W'_0$, and the stabilizer of $\Upsilon_0$ is $W_0$, so 
$\Sigma_0$ is completely determined by an element in $W/W_0$. Proceeding inductively,
we see \cite{CC}:
\begin{prop}
The set of all the combinatorial galleries of type $t_{\g_{\lam}}$ and of source $F_f$, $\G(\g_\lam)$, 
is in bijection with the quotient
$$
W\times_{W_0} W'_1\times_{W_1} \cdots \times_{W_{p-1} }W'_p/{W_p},
$$
of the group $W\times W'_1\times \cdots \times W'_p$ by the subgroup 
${W_0} \times W_1\times  \cdots \times{W_p}$
under the action defined by $(w'_0 w_0, w_0^{-1} w'_1 w_1,..., w_{p-1}^{-1} w'_p w_p)$.
\end{prop} 

The image of a gallery $\gamma$ will be denoted by $\gamma=[\de_0,\de_1,...,\de_p]$. 
Because of the definition of the quotient we can (and will) assume that $\de_j\in W_j'$ is 
{\it the unique representative of minimal length of its class in} $W'_j/W_j$. 
We will use freely both of the notation for a gallery. 
\begin{exam}\label{taum}\rm
The gallery $\g_\lam$ can be written as $\g_\lam = [1,\tau^m_1,...,\tau^m_p]$, where each 
$\tau^m_j\in W'_j$  is the {\it unique minimal representative of the largest class} 
(in the induced Bruhat ordering) in  $W'_j/W_j$. The product 
$\tau^m_1\cdots\tau^m_p\in\Waff$ is a reduced decomposition of the element 
that sends $\De_f$ to $proj_{F_\lam}(\De_f)$. 
All galleries of shape $[\g_0,\tau^m_1,...,\tau^m_j]$ for 
$\g_0\in W$ (of minimal length modulo $W_\lam$) are minimal.
\end{exam}

Let $\de = [\de_0,\de_1,...,\de_p] = (F_f\subset \Sigma_0\supset {\Sigma_1}'\subset\cdots 
\supset{\Sigma}'_p \subset {\Sigma}_{p}\supset F_{\nu})\in\G(\g_\lam)$. The gallery $\de$ 
is called {\it folded around  the small face $\Sigma'_j$} if $\de_j \ne \tau^m_j$. 
Such a gallery will be obtained from $[\de_0,\tau^m_1,...,\tau^m_p]$ by applying at the places 
$j$, where $\de_j\ne\tau^m_j$, some affine reflections with respect to the affine hyperplanes containing
the small faces. I.e., associated to $\de$ we have the following sequence of galleries in $ \G(\g_\lam)$:
\begin{equation}\label{sequenceoffoldings}
\begin{array}{crclrcl}

 \g_0=[\de_0,\tau^m_1,\tau^m_2,...,\tau^m_p],&\g_1&=&[\de_0,\de_1,\tau^m_2,...,\tau^m_p],      
 &\g_2&=&[\de_0,\de_1,\de_2,\tau^m_3,...,\tau^m_p],\\
\ldots                                                                 &\g_{p-1}&=&[\de_0,\de_1,...,\de_{p-1},\tau^m_p],&\g_{p}&=&\de.
\end{array}
\end{equation}
Here $\g_1$ is obtained from the minimal gallery $\g_0=\de$ by a folding around 
${\Sigma}'_1$, $\g_2$ is obtained by a folding of $\g_1$ around ${\Sigma}'_2$, etc.

\begin{dfn}\rm
The gallery $\de\in \G(\g_\lam)$ is called {\it positively folded} at $\Sigma'_j$ if for all the affine hyperplanes 
$\ch$ involved in the sequence of reflections for the folding of $\g_{j-1}$ around $\Sigma_j'$ to get $\g_j$, 
the image of the reflection is ``separated'' by $\ch$ from the anti--dominant
chamber $\Lc_{-\infty}$ at $-\infty$. We say that the gallery is {\it positively folded} if all foldings
are positive.
\end{dfn}
Here ``{\it separated\/}'' means that there {\it exists a representative $\Ls$} 
of the class $\Lcinf$ such that the {\it image of the reflection and $\Ls$ are separated} by the reflection hyperplane.
\begin{exam}\rm
The minimal galleries are automatically positively folded.
\end{exam}

Let $\Gamma^+(\g_\lam)$ denote the subset of $\Gamma (\g_\lam)$ of all {\it positively 
folded combinatorial galleries of type $t_{\gamma_\lam}$}, and for a co-character $\nu$, we denote 
$\Gamma^+(\g_\lam,\nu)$ the subset of galleries $\de$ in $\Gamma^+(\g_\lam)$ ending in $\nu$, 
i.e., $\de= (F_f\subset \Sigma_0\supset \ldots \subset \Sigma_{p}\supset F_{\nu})$. 
\begin{exam}\label{exampleSl3}
\rm
Consider the group $G=SL_3(\bc)$, let $\Saff=\{s_0,s_1,s_2\}$ be indexed such that
the Weyl group $W$ of $G$ is generated by $\{s_1,s_2\}$. A minimal gallery joining
the origin with the highest root $\beta$ is $\gamma_\beta=[1,s_0]$. The elements in the set 
$\Gamma^+(\g_\beta)$  of all positively folded galleries of the same type and source $0$ are:
$$
I=[1,s_0],\,I\hskip -2pt I=[s_1,s_0],\,I\hskip -2pt I\hskip -2pt I=[s_2,s_0],\,I\hskip -2pt V= 
[s_2s_1,s_0],\, V=[s_1s_2,s_0],\,V\hskip -2pt I= [s_2s_1s_2,s_0]
$$
which are the minimal galleries in this set, and
$$
V\hskip -2pt I\hskip -2pt I=[s_1s_2,1],\,V\hskip -2pt I\hskip -2pt I\hskip -2pt I=
[s_2s_1,1],\,I\hskip -2pt X=[s_2s_1s_2,1].
$$
\hskip 3cm
\beginpicture \footnotesize \setlinear
\setcoordinatesystem units  <6mm, 6mm>
\setplotarea x from -4.5 to 4.5, y from -4.0 to 4.0
\accountingoff 
\setdots<5pt>
\plot -5  1.73   -4 3.46 /
\plot -5 -1.73  -2 3.46 /
\plot -4 -3.46   0 3.46 /
\plot -2 -3.46   1  1.73 / 
\plot  0 -3.46   1 -1.73 /  
\plot -4. -3.46   -5 -1.73 /
\plot -2 -3.46   -5 1.73 /
\plot 0 -3.46   -4 3.46 /
\plot 1 -1.73   -2 3.46 /
\plot 1  1.73    0 3.46 /
\plot -5.00 0.00    1.00 0 /
\plot -5.00 1.73    1.00 1.73 /
\plot -5.00 3.46    1.00 3.46 /
\plot -5.00 -1.73   1.00 -1.73 /
\plot -5.00 -3.46   1.00 -3.46 /
\plot -5.00 -3.46  -5.00 3.46 /
\plot  1.00 -3.46   1.00 3.47 /
\setdots<0pt>
\plot -3.95 -0.05    -2.05 -0.05 /
\plot -4 -0.05    -3 -1.73 /
\plot -3 -1.65       -2.1 -0.05  /
\plot -4 -0.05    -5 -1.73 /
\plot -3.1 -1.73      -5 -1.73  /

\plot -3.95 0.05    -2.05  0.05 /
\plot -4 0.05    -3  1.73 /
\plot -3  1.65       -2.1  0.05  /
\plot -4  0.05    -5  1.73 /
\plot -3.1  1.73      -5  1.73  /

\plot -0.05 -0.05    -1.95 -0.05 /
\plot -0 -0.05    -1 -1.73 /
\plot -1 -1.65       -1.9 -0.05  /
\plot -0 -0.05     1 -1.73 /
\plot -0.9 -1.73       1 -1.73  /

\plot -0.05  0.05    -1.95  0.05 /
\plot -0  0.05    -1  1.73 /
\plot -1  1.65       -1.9  0.05  /
\plot -0  0.05     1  1.73 /
\plot -0.9  1.73       1  1.73  /

\plot -2.8 -1.73       -2 -0.25  /
\plot -2.8 -1.73       -1.2 -1.73   /
\plot -2 -0.25       -1.2 -1.65   /
\plot -2.8 -1.73       -2 -3.30   /
\plot -2 -3.30      -1.2 -1.73   /

\plot -2.8 1.73       -2  0.25  /
\plot -2.8 1.73       -1.2  1.73   /
\plot -2 0.25       -1.2  1.65   /
\plot -2.8 1.73       -2  3.30   /
\plot -2 3.30      -1.2  1.73   /

\put{\hbox{$\epsilon_1-\epsilon_2$}} at 1.8 1.2
\put{\hbox{$\epsilon_2-\epsilon_3$}} at -5.8 1.2
\put{\hbox{$\epsilon_1-\epsilon_3$}} at -.8 3.2
\put{\hbox{$I$}} at -2 2.1
\put{\hbox{$I\hskip -2pt I$}} at -4 1.3
\put{\hbox{$I\hskip -2pt I\hskip -2pt I$}} at 0 1.3
\put{\hbox{$I\hskip -2pt V$}} at 0 -1.3
\put{\hbox{$V$}} at -4 -1.3
\put{\hbox{$V\hskip -2pt I$}} at -2 -2.1
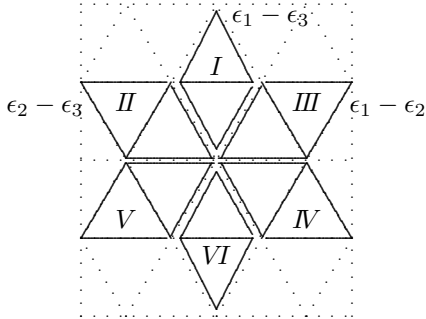
\put{\hbox{\rm The minimal galleries}} at -2 -4

\hskip 200pt
\setdots<5pt>
\plot -5.00  1.73   -4.00 3.46 /
\plot -5.00 -1.73  -2.00 3.46 /
\plot -4.00 -3.46   0.00 3.46 /
\plot -2.00 -3.46   1  1.73 / 
\plot  0.00 -3.46   1 -1.73 /  
\plot -4.00 -3.46   -5.00 -1.73 /
\plot -2.00 -3.46   -5.00 1.73 /
\plot 0.00 -3.46   -4.00 3.46 /
\plot 1.00 -1.73   -2.00 3.46 /
\plot 1.00  1.73    0.00 3.46 /
\plot -5.00 0.00    1.00 0 /
\plot -5.00 1.73    1.00 1.73 /
\plot -5.00 3.46    1.00 3.46 /
\plot -5.00 -1.73   1.00 -1.73 /
\plot -5.00 -3.46   1.00 -3.46 /
\plot -5.00 -3.46  -5.00 3.46 /
\plot  1.00 -3.46   1.00 3.47 /
\setdots<0pt>
\plot -3.95 -0.05    -0.05 -0.05 /
\plot -3.90 -0.05    -3 -1.65 /
\plot -3 -1.65       -2.1 -0.05  /
\plot -1.90 -0.05    -1 -1.65 /
\plot -1 -1.65       -0.1 -0.05  /
\plot -2.8 -1.65       -2 -0.25  /
\plot -2.8 -1.65       -1.2 -1.65   /
\plot -2 -0.25       -1.2 -1.65   /
\setdashes
\plot -3.95 -0.1    -0.05 -0.1 /
\plot -3.80 -0.1    -3 -1.50 /
\plot -3 -1.50       -2.2 -0.05  /
\plot -1.80 -0.1    -1 -1.50 /
\plot -1 -1.50       -0.2 -0.05  /
\plot -2.7 -1.65       -2 -0.35  /
\plot -2.8 -1.55       -1.2 -1.55   /
\plot -2 -0.4            -1.35 -1.55   /
\put{\hbox{$0$}} at -2 .3
\put{\hbox{$V\hskip -2pt I\hskip -2pt I\hskip -2pt I$}} at -1 -.45
\put{\hbox{$V\hskip -2pt I\hskip -2pt I$}} at -3 -.45
\put{\hbox{$I\hskip -3pt X$}} at -2 -1.25
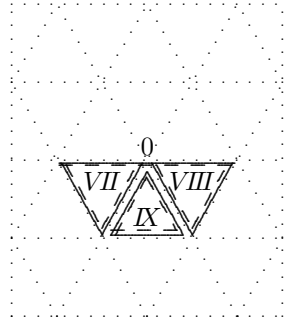
\put{\hbox{The folded galleries}} at -2 -4
\endpicture

\end{exam}

Note: for each root we have exactly one such gallery ending in this weight (the minimal galleries), 
and for the zero weight we have three (the non-minimal galleries), which is one too much if one has 
expected to get the character of the adjoint representation of $SL_3(\bc)$. We will see in the next 
section how to refine the choice.


\section{Dimension of galleries}\label{dimensiongallery}
We will see that the notion of the dimension of a gallery is related to the dimension of the fibre
of the retraction $\rinf$ with center at $-\infty$ discussed later. 
Fix a positively folded gallery  $\g\in \Gamma^+(\g_\lam)$ of type $t_{\gamma_\lam}$:
$$
\g = (F_f = \G'_0\subset \Gamma_0\supset \Gamma'_1\subset\cdots\subset\G_{j-1}\supset
\G'_j\subset\G_j\supset\cdots\supset\Gamma'_p \subset \Gamma_{p}\supset F_{\nu}).
$$  
For $j = 0,...,p$, let $\sh_j$ be the set of all the affine hyperplanes $\ch$ in $\cha$ such 
that $\G'_j\subset\ch$ and $\G_j\not\subset \ch$. 
We say that an affine hyperplane $\ch$
is a {\it load-bearing wall} for $\g$ {\it at} $\Gamma_j$ if  
$\ch\in\sh_j$ and $\ch$ separates $\Gamma_{j}$ from $\Lc_{-\infty}$.  
Note: for a positively folded gallery all folding hyperplanes are load-bearing walls by definition.

\begin{dfn}\rm
The {\it dimension} of the gallery $\g\in\Gamma^+(\g_\lam)$ is the number of pairs
$(\ch,\Gamma_j)$ such that $\ch$ is a load-bearing wall
for $\g$ at $\Gamma_j$:
$$
\dim \g=\sharp\{(\ch,\Gamma_j)\mid\ch\hbox{\rm\ is a load-bearing wall for $\g$ at \ }
\Gamma_j\}.
$$
\end{dfn}

\begin{exam}\label{exampleextremalgallery}\rm
Let $\lam$ be a dominant co-character and let $\g_\lam$ be a minimal gallery
joining $0$ with $\lam$. If $\beta$ is a positive root such that $\la\lam,\be\ra > 0$, then the affine hyperplanes
$\ch_{\beta,n}$ for all $0\le n<\la\lam,\be\ra$ are load-bearing walls. It follows that 
$$
\dim \gamma_\lam=\sum_{\beta\in\Phi^+}\la\lam,\be\ra=\la\lam,2\rho\ra=\la2\lam,\rho\ra.
$$ 
Next, let $\nu=\tau(\lam)$ be a Weyl group conjugate of $\lam$ (as usual, we identify  
$\tau\in W/W_\lam$ with its representative of minimal length), and let $\gamma_\tau$
be the gallery obtained from $\g_\lam$ by applying simultaneously $\tau$ to all
faces of $\g_\lam$. Since $0 = F_f$ is a fixed point for the action of $W$ and $\nu=\tau(\lam)$,
this is a gallery joining $0$ with $\nu$. By the minimality of $\gamma_\lam$, this
gallery is also minimal and has no foldings. So, in particular, $\gamma_\tau\in \Gamma^+(\g_\lam)$. 
One shows easily by decreasing induction on $\la\nu,\rho\ra$ that 
$$
\dim\g_\tau=\la\lam+\nu,\rho\ra.
$$
\end{exam}
Our next aim is to provide such a formula for all $\g\in\Gamma^+(\g_\lam)$. Recall
that $\Gamma^+(\g_\lam,\nu)$ is the set of all galleries $\g$ in $\Gamma^+(\g_\lam)$
ending in $\nu$, i.e., $\g= (F_f\subset \Gamma_0\supset \cdots
\subset \Gamma_{p+1}\supset F_{\nu})$.

\begin{prop}\label{dim-inequality}
If $\gamma\in \Gamma^+(\g_\lam,\nu)$, then $\dim\g\le \la\lam+\nu,\rho\ra$.
\end{prop}

\begin{dfn}\rm
A positively folded gallery $\gamma\in \Gamma^+(\g_\lam,\nu)$ is called an 
{\it LS-gallery of type} $\lam$ if  $\dim\g = \la\lam+\nu,\rho\ra$.
\end{dfn}
The proofs of the propositions~\ref{dim-inequality} and \ref{characterformula} (see below) 
will be given in the next section, where we also provide a combinatorial
characterization of the LS-galleries. Denote $\Gamma_{LS}^+(\g_\lam,\nu)$ the set
of all galleries in $\Gamma^+(\g_\lam,\nu)$ that are LS-galleries.
As a consequence of the description we will get the following combinatorial
character formula, which 
is the equivalent in terms of galleries of the path character 
formula in \cite{Lit1,Lit2}: Let $V(\lam)$  be the irreducible complex representation of highest
weight $\lam$ for the semisimple algebraic group $G^\vee$ (the Langland's dual group
of $G$), and denote $\chara V(\lam)$ its character:
\begin{prop}\label{characterformula}
$\chara V(\lam)=\sum_{\nu\in \Xvee} 
\sharp\Gamma_{LS}^+(\g_\lam,\nu)\exp(\nu)$.
\end{prop}


\section{Root Operators}\label{rootoperators}
The aim of this section is to prove the dimension formula and to give a combinatorial
characterization of the LS-galleries. We define now ``folding" operators $f_\al,e_\al,\ce_\al $
(for all simple roots) on the set of all combinatorial galleries $\Gamma(\g_\lam)$ of a fixed type.
Let $\lam$ be a dominant co-character and fix a minimal gallery $\g_\lam$  joining
the origin with $\lam$.

Let $\al$ be a simple root and let $\nu\preceq \lam$ be a co-character. Fix a combinatorial gallery
$\g\in\Gamma(\g_\lam,\nu)$ of type $t_{\gamma_\lam}$, say
$$
\g = [\g_0,\g_1,...,\g_p] = (F_f=\G_0'\subset \Gamma_0\supset \Gamma'_1\subset\cdots
\subset \Gamma_{p}\supset \G'_{p+1}=F_{\nu}).
$$
Let $m\in\bz$ be minimal such that one of the small faces $\Gamma_k'$
is contained in the hyperplane $\ch_{\al, m}$, note that $m\le 0$. The
operators $e_\al$ and $\ce_\al$ are different, but the conditions for
the operators to be defined are not exclusive, so it might well be that both are
defined for a given gallery. 
\begin{itemize}
\item [{I)}] Suppose that $m\le -1$. Let $k$ be minimal such that $\Gamma_k'\subset \ch_{\al, m}$,
and fix $0\le j\le k$ maximal such that the small face $\Gamma_j'\subset \ch_{\al, m+1}$ is contained
in the hyperplane $\ch_{\al, m+1}$.
\end{itemize}

\hskip 1cm
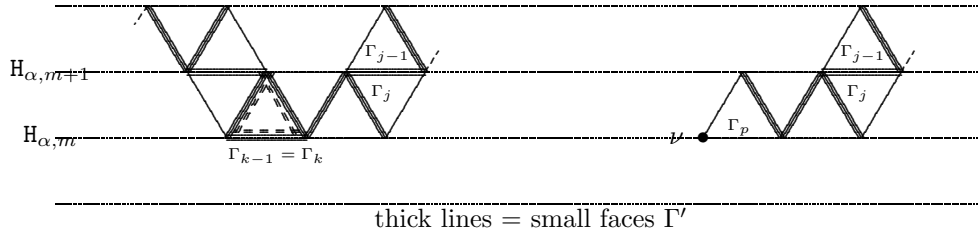
\beginpicture \footnotesize \setlinear
\setcoordinatesystem units  <100pt, 100pt>
\setplotarea x from -2.0 to 1.5, y from -.75 to .26
\setdots<1pt>
\axis bottom shiftedto y=0 /
\axis bottom shiftedto y=.25 /
\axis bottom shiftedto y=-.25 /
\axis bottom shiftedto y=-.5 / 
\setsolid                                           
\plot -.6 0  -1.5   0 /
\plot  -.61  .0  -.76   .25 / 
\plot  -.59  .0  -.74    .25 /
\plot  -1.51  .0  -1.66   .25 / 
\plot  -1.49  .0  -1.64    .25 /
\setdashes<2pt>
\plot  -.6  0   -.55   .08 / 
\plot  -.75  .25   -.65   .25 /
\plot  -1.65  .25   -1.7   .18 / 
\plot  -1.5  0   -1.6   0 /
\setsolid
\plot  -.6   .01  -.9    .01 /
\plot  -.6  -.01  -.9   -.01 / 
\plot  -.91   0  -.76    -.25 /
\plot  -.89   0  -.74   -.25 / 
\plot  -.91   0  -1.06    -.25 /
\plot  -.89   0  -1.04   -.25 / 
\plot  -1.5   .01  -1.2    .01 /
\plot  -1.5  -.01  -1.2   -.01 /
\plot  -1.21   .0  -1.06    -.25 /
\plot  -1.19  -.0  -1.04    -.25 /
\plot  -1.21   .0  -1.36    -.25 /
\plot  -1.19  -.0  -1.34   -.25  /
\plot  -1.51   0  -1.36    .25 /
\plot  -1.49   0  -1.34    .25 /
\plot  -1.35   -.24  -1.05    -.24 /
\plot  -1.35   -.26  -1.05    -.26 /
\plot -.75  -.25  -1.35  -.25 /
\plot -.75   .25  -1.05  -.25 /
\plot -.75   .25 -.6 0 /
\plot -.9      0 -.75   -.25  /
\plot -.6  0 -.75   -.25  /
\plot -1.05  -.25  -1.35  .25  /
\plot -1.35  -.25  -1.65  .25  /
\plot -1.35  -.25  -1.2  0  /
\plot -1.5  0  -1.35  .25  /
\plot -1.65  .25  -1.35  .25  /
\setdashes <3pt>
\plot -1.22   0  -1.08  -.25  /
\plot -1.23   0  -1.09  -.25  /
\plot -1.18   0  -1.33  -.25  /
\plot -1.17   0  -1.32  -.25  /
\plot  -1.076  -.22   -1.32  -.22  /
\plot  -1.076  -.23   -1.32  -.23  /
\setsolid
\put{\hbox{\bf $\ch_{\al,m+1}$ }} at   -2  0
\put{\hbox{\bf $\ch_{\al,m}$ }} at   -2  -.25
\put{\tiny{\bf $\Gamma_{j-1}$ }} at   -.73  .07
\put{\tiny{\bf $\Gamma_{j}$ }} at   -.75  -.08
\put{\tiny{\bf $\Gamma_{k-1}=\Gamma_{k}$ }} at   -1.15  -.32
\put{\hbox{thick lines = small faces $\Gamma'$}} at  -0.2  -.55
\plot  1.2  0  .6   0 /
\plot  1.06  .25   1.21   0 / 
\plot  1.04  .25   1.19   0 /
\setdashes<2pt>
\plot  1.2  0   1.25   .08 / 
\plot  1.05  .25   1.15   .25 /
\setdots<0pt>
\plot  1.2  -.01  .9   -.01 / 
\plot  1.2   .01  .9    .01 /
\plot   .89   0   1.04   -.25 / 
\plot   .91   0   1.06     -.25 /
\plot   .89   0   .74   -.25 / 
\plot   .91   0   .76     -.25 /
\plot   .59   0   .74   -.25 / 
\plot   .61   0   .76     -.25 /
\plot  1.05  .25  .75  -.25 /
\plot 1.05  .25  1.2  0 /
\plot 1.2  0 1.05 -.25 /
\plot .9  0 1.05 -.25 /
\plot .45  -.25 1.05 -.25 /
\plot .45 -.25 .6 0 /
\plot .75 -.25 .6 0 /
\put{\tiny{\bf $\Gamma_{j-1}$ }} at   1.07  .07
\put{\tiny{\bf $\Gamma_{j}$ }} at   1.05  -.08
\put{\tiny{\bf $\Gamma_{p}$ }} at   .6  -.21
\put{$\bullet$} at   .45  -.25
\put{$\nu$} at   .35  -.25
\endpicture

\vskip-0.5cm
\begin{dfn}\label{defealpha}\rm
If $\g = [\g_0,\g_1,...,\g_p]$ satisfies I), then let $e_\al \g$ be the gallery defined by:
$$
e_\al \g = [\de_0,\de_1,...,\de_p] = (F_f\subset \De_0\supset
\De'_1\subset\cdots\supset \De'_p \subset
\De_{p}\supset  F_{\check\nu})
$$
$$
\hbox{\rm where }\quad
\De_i = \left\{
\begin{array}{ll}
\Gamma_i & \hbox{ for } i\le j-1,\\
s_{\al,m+1}(\Gamma_i) & \hbox{ for } j \leq i\leq k-1,\\
t_{\avee}(\Gamma_i) &
\hbox{ for } i\geq k.
\end{array}\right\}\
\hbox{\rm and }\  t_{\avee}=\hbox{\rm translation by}\ \avee.
$$
\end{dfn}
We define a partial inverse operator $f_\al$ to the operator $e_\al$.
\begin{itemize}
\item [{II)}] Suppose that $\la\nu,\al\ra-m\ge 1$. Let $j$ be maximal such that
 $\Gamma_j'\subset \ch_{\al, m}$ and fix $j\le k\le p+1$ minimal such that the small face
 $\Gamma_k'$ is contained in $\ch_{\al, m+1}$.
\end{itemize}

\hskip 1cm
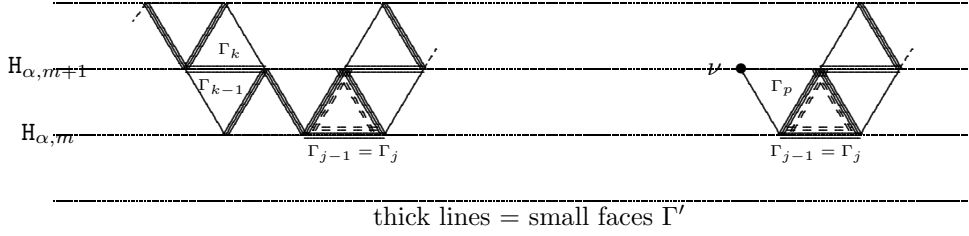
\hbox{
\beginpicture \footnotesize \setlinear
\setcoordinatesystem units  <100pt, 100pt>
\setplotarea x from -2.0 to 1.5, y from -.75 to .26
\setdots<1pt>
\axis bottom shiftedto y=0 /
\axis bottom shiftedto y=.25 /
\axis bottom shiftedto y=-.25 /
\axis bottom shiftedto y=-.5 / 
\setdots<0pt>                                            
\plot -.6 0  -1.5   0 /
\plot  -.61  .0  -.76   .25 / 
\plot  -.59  .0  -.74    .25 /
\plot  -1.51  .0  -1.66   .25 / 
\plot  -1.49  .0  -1.64    .25 /
\setdashes<2pt>
\plot  -.6  0   -.55   .08 / 
\plot  -.75  .25   -.65   .25 /
\plot  -1.65  .25   -1.7   .18 / 
\plot  -1.5  0   -1.6   0 /
\setdots<0pt>
\plot  -.6   .01  -.9    .01 /
\plot  -.6  -.01  -.9   -.01 / 
\plot  -.91   0  -.76    -.25 /
\plot  -.89   0  -.74   -.25 / 
\plot  -.91   0  -1.06    -.25 /
\plot  -.89   0  -1.04   -.25 / 
\plot  -1.5   .01  -1.2    .01 /
\plot  -1.5  -.01  -1.2   -.01 /
\plot  -1.21   .0  -1.06    -.25 /
\plot  -1.19  -.0  -1.04    -.25 /
\plot  -1.21   .0  -1.36    -.25 /
\plot  -1.19  -.0  -1.34   -.25  /
\plot  -1.51   0  -1.36    .25 /
\plot  -1.49   0  -1.34    .25 /
\plot -.75  -.25  -1.35  -.25 /
\plot -.75   .25  -1.05  -.25 /
\plot -.75   .25 -.6 0 /
\plot -.9      0 -.75   -.25  /
\plot -.6  0 -.75   -.25  /
\plot -1.05  -.25  -1.35  .25  /
\plot -1.35  -.25  -1.65  .25  /
\plot -1.35  -.25  -1.2  0  /
\plot -1.5  0  -1.35  .25  /
\plot -1.65  .25  -1.35  .25  /
\setdashes <3pt>
\plot -.92   0  -.78  -.25  /
\plot -0.93   0  -0.79  -.25  /
\plot -0.88   0  -1.03  -.25  /
\plot -0.87   0  -1.02  -.25  /
\plot  -0.776  -.22   -1.02  -.22  /
\plot  -0.776  -.23   -1.02  -.23  /
\setdots<0pt>
\plot  -1.05   -.24  -0.75    -.24 /
\plot  -1.05   -.26  -0.75    -.26 /
\put{\hbox{\bf $\ch_{\al,m+1}$ }} at   -2  0
\put{\hbox{\bf $\ch_{\al,m}$ }} at   -2  -.25
\put{\tiny{\bf $\Gamma_{j-1}=\Gamma_{j}$ }} at   -0.85  -.32
\put{\tiny{\bf $\Gamma_{k-1}$ }} at   -1.35  -.07
\put{\tiny{\bf $\Gamma_{k}$ }} at   -1.32  .07
\put{\hbox{thick lines = small faces $\Gamma'$}} at  -0.2  -.55
\plot  1.2  0  .6   0 /
\plot  1.06  .25   1.21   0 / 
\plot  1.04  .25   1.19   0 /
\setdashes<2pt>
\plot  1.2  0   1.25   .08 / 
\plot  1.05  .25   1.15   .25 /
\setdots<0pt>
\plot  1.2  -.01  .9   -.01 / 
\plot  1.2   .01  .9    .01 /
\plot   .89   0   1.04   -.25 / 
\plot   .91   0   1.06     -.25 /
\plot   .89   0   .74   -.25 / 
\plot   .91   0   .76     -.25 /
\plot  1.05  .25  .75  -.25 /
\plot 1.05  .25  1.2  0 /
\plot 1.2  0 1.05 -.25 /
\plot .9  0 1.05 -.25 /
\plot .75  -.25 1.05 -.25 /
\plot .75  -.24 1.05 -.24 /
\plot .75  -.26 1.05 -.26 /
\setdashes <3pt>
\plot  .92   0       .78  -.25  /
\plot  0.93   0   0.79  -.25  /
\plot  0.88   0   1.03  -.25  /
\plot  0.87   0   1.02  -.25  /
\plot   0.776  -.22    1.02  -.22  /
\plot   0.776  -.23    1.02  -.23  /
\setdots<0pt>
\plot .75 -.25 .6 0 /
\put{\tiny{\bf $\Gamma_{j-1}=\Gamma_{j}$ }} at   .9  -.32
\put{\tiny{\bf $\Gamma_{p}$ }} at   .77  -.07
\put{$\bullet$} at   .6  -0
\put{$\nu$} at   .5  0
\endpicture
}
\vskip-0.5cm

\begin{dfn}\rm
If $\g$ satisfies II), then let $f_\al \g$ be the gallery defined by:
$$
f_\al \g = [\g_0,\g_1,...,\g_p] = (F_f\subset \De_0\supset
\De'_1\subset\cdots\supset \De'_p \subset
\De_{p}\supset  F_{\check\nu})
$$
$$
\hbox{\rm where}\quad
\De_i = \left\{
\begin{array}{ll}
\Gamma_i & \hbox{ for } i< j,\\ s_{\al,m}(\Gamma_i) & \hbox{ for } j \leq i< k,
\\ t_{-\avee}(\Gamma_i) &
\hbox{ for } i\geq k.
\end{array}\right.
$$
\end{dfn}

\begin{itemize}
\item[{III)}] Assume that $\gamma$ crosses $\ch_{\al, m}$.
Fix $j$ minimal with this property, i.e., $\Gamma_j'\subset \ch_{\al, m}$ and
$\ch_{\al, m}$ separates the faces $\Gamma_i$ from $\Lc_{-\infty}$ for $i<j$
(in the large sense, i.e., $\Gamma_i\subset\ch_{\al, m}^+$) but not the face
$\Gamma_{j}$ (so $\ch_{\al, m}\in \sh_j$, the set of all the affine hyperplanes
$\ch$ such that $\G'_j\subset\ch$ and $\G_j\not\subset \ch$). Fix $k>j$ minimal such that
$\Gamma_k'\subset \ch_{\al, m}$.
\end{itemize}

\hskip 1cm
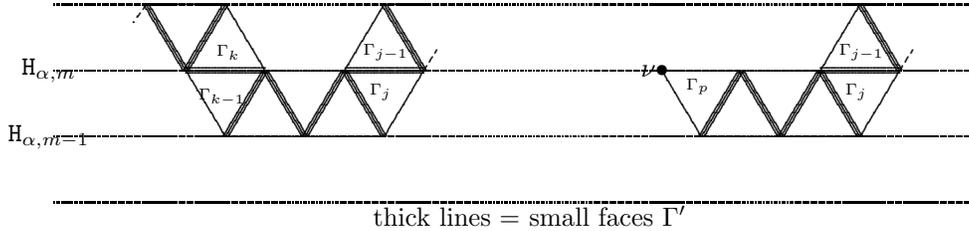
\hbox{
\beginpicture \footnotesize \setlinear
\setcoordinatesystem units  <100pt, 100pt>
\setplotarea x from -2.0 to 1.5, y from -.75 to .26
\setdots<1pt>
\axis bottom shiftedto y=0 /
\axis bottom shiftedto y=.25 /
\axis bottom shiftedto y=-.25 /
\axis bottom shiftedto y=-.5 / 
\setdots<0pt>                                            
\plot -.6 0  -1.5   0 /
\plot  -.61  .0  -.76   .25 / 
\plot  -.59  .0  -.74    .25 /
\plot  -1.51  .0  -1.66   .25 / 
\plot  -1.49  .0  -1.64    .25 /
\setdashes<2pt>
\plot  -.6  0   -.55   .08 / 
\plot  -.75  .25   -.65   .25 /
\plot  -1.65  .25   -1.7   .18 / 
\plot  -1.5  0   -1.6   0 /
\setdots<0pt>
\plot  -.6   .01  -.9    .01 /
\plot  -.6  -.01  -.9   -.01 / 
\plot  -.91   0  -.76    -.25 /
\plot  -.89   0  -.74   -.25 / 
\plot  -.91   0  -1.06    -.25 /
\plot  -.89   0  -1.04   -.25 / 
\plot  -1.5   .01  -1.2    .01 /
\plot  -1.5  -.01  -1.2   -.01 /
\plot  -1.21   .0  -1.06    -.25 /
\plot  -1.19  -.0  -1.04    -.25 /
\plot  -1.21   .0  -1.36    -.25 /
\plot  -1.19  -.0  -1.34   -.25  /
\plot  -1.51   0  -1.36    .25 /
\plot  -1.49   0  -1.34    .25 /
\plot -.75  -.25  -1.35  -.25 /
\plot -.75   .25  -1.05  -.25 /
\plot -.75   .25 -.6 0 /
\plot -.9      0 -.75   -.25  /
\plot -.6  0 -.75   -.25  /
\plot -1.05  -.25  -1.35  .25  /
\plot -1.35  -.25  -1.65  .25  /
\plot -1.35  -.25  -1.2  0  /
\plot -1.5  0  -1.35  .25  /
\plot -1.65  .25  -1.35  .25  /
\put{\hbox{\bf $\ch_{\al,m}$ }} at   -2  0
\put{\hbox{\bf $\ch_{\al,m-1}$ }} at   -2  -.25
\put{\tiny{\bf $\Gamma_{j-1}$ }} at   -.73  .07
\put{\tiny{\bf $\Gamma_{j}$ }} at   -.75  -.08
\put{\tiny{\bf $\Gamma_{k-1}$ }} at   -1.35  -.1
\put{\tiny{$\Gamma_{k}$}} at   -1.34  .07
\plot  1.2  0  .3   0 /
\plot  .3   0 .45 -.25 /
\plot  1.06  .25   1.21   0 / 
\plot  1.04  .25   1.19   0 /
\setdashes<2pt>
\plot  1.2  0   1.25   .08 / 
\plot  1.05  .25   1.15   .25 /
\setdots<0pt>
\plot  1.2  -.01  .9   -.01 / 
\plot  1.2   .01  .9    .01 /
\plot   .89   0   1.04   -.25 / 
\plot   .91   0   1.06     -.25 /
\plot   .89   0   .74   -.25 / 
\plot   .91   0   .76     -.25 /
\plot   .59   0   .74   -.25 / 
\plot   .61   0   .76     -.25 /
\plot   .59   0   .44   -.25 / 
\plot   .61   0   .46     -.25 /
\plot  1.05  .25  .75  -.25 /
\plot 1.05  .25  1.2  0 /
\plot 1.2  0 1.05 -.25 /
\plot .9  0 1.05 -.25 /
\plot .45  -.25 1.05 -.25 /
\plot .45 -.25 .6 0 /
\plot .75 -.25 .6 0 /
\put{\tiny{\bf $\Gamma_{j-1}$ }} at   1.07  .07
\put{\tiny{\bf $\Gamma_{j}$ }} at   1.05  -.08
\put{\tiny{\bf $\Gamma_{p}$ }} at   .45  -.06
\put{$\bullet$} at   .3  0
\put{$\nu$} at   .25  0
\put{\hbox{thick lines = small faces $\Gamma'$}} at  -0.2  -.55
\endpicture
}
\vskip-0.5cm

\begin{dfn}\label{defetildealpha}\rm
If $\g = [\g_0,\g_1,...,\g_p]$ satisfies III), then let $\widetilde{e}_\al \g$ be the gallery defined by:
$$
\widetilde{e}_\al \g = [\de_0,\de_1,...,\de_p] = (F_f\subset \De_0\supset
\De'_1\subset\cdots\supset \De'_p \subset
\De_{p}\supset  F_{\check\nu})
$$ 
$$
\hbox{\rm where}\quad \De_i = \left\{
\begin{array}{ll}
\Gamma_i & \hbox{ for } i\le j-1  \hbox{ and } i\geq k,\\
s_{\al,m}(\Gamma_i) & \hbox{ for } j \leq i< k.
\end{array}\right.
$$
\end{dfn}
It follows immediately from the definition of the operators that if one of the galleries $e_\al\g$,$f_\al\g$ or $\widetilde{e}_\al\g$ is defined,
then it is again a gallery of type $\g_\lam$.

The proof of the following simple lemma is left to the reader.
Let $\g\in\Gamma(\g_\lam,\nu)$ be a combinatorial gallery
of type $t_{\g_\lam}$ ending in $\nu$, let $\al$ be a simple root and suppose that $m$ is
minimal such that one of the small faces $\Gamma_k'\subset \ch_{\al, m}$
is contained in the hyperplane $\ch_{\al, m}$
\begin{lem}\label{simplelemmaA}
\begin{itemize}
\item[({\it i\/})] The gallery $e_\al\g$ is not defined if and only if $ m=0$,
and if $e_\al\g$ is defined, then $ e_\al\g\in\Gamma(\g_\lam,\nu+\avee)$.
\item[({\it ii\/})] The gallery $f_\al\g$ is not defined if and only if $ m=\la\nu,\al\ra$,
and if $f_\al\g$ is defined, then $ f_\al\g\in\Gamma(\g_\lam,\nu-\avee)$.
\item[({\it iii\/})] If $e_\al\g$ is defined, then $f_\al(e_\al\g)$ is
defined and equal to $\g$. Further, $m+1$ is minimal such that
a small face of the gallery $e_\al\g$ is contained in $\ch_{\al,m+1}^+$.
\item[({\it iv\/})] If $f_\al\g$ is defined, then $e_\al(f_\al\g)$ is
defined and equal to $\g$. Further, $m-1$ is maximal such that
a small face of the gallery $f_\al\g$ is contained in $\ch_{\al,m-1}^+$.
\item[({\it v\/})] Let $p$ be maximal such that $f_\al^p\g$  is defined
and let $q$ be maximal such that $e_\al^q\g$ is defined, then
$p-q=\la\nu,\al\ra$.
\end{itemize}
\end{lem}

Let $\G(\g_\lam,{\rm dom})$ be the set of all combinatorial galleries $\de$
in $\G(\g_\lam)$ such that none of the small faces is contained in one
of the hyperplanes $\ch_{\al,-1}$, or, in other words, $e_\al\de$ is not
defined for all simple roots. Set $\chara \G(\g_\lam)=\sum \exp(e(\g))$,
the sum over all $\g\in \G(\g_\lam)$. 

\begin{coro}\label{characterdecompo}
$\chara \G(\g_\lam) = \sum \chara V(e(\g))$, the sum running over all
$\g\in \G(\g_\lam,{\rm dom})$.
\end{coro}

\proof
The character formula can be proved using the same arguments 
as in \cite{Lit5}, so we will only give a sketch of the proof.
It follows easily from Lemma~\ref{simplelemmaA} that $\chara \G(\g_\lam)$
is stable under the Weyl group. Using Weyl's character formula, the above is equivalent to
$$
(\sum_{w\in W}\sgn(w)e^{w(\rho)})\chara \G(\g_\lam)
=\sum_{\scriptstyle \g\in \G(\g_\lam,{\rm dom})}
(\sum_{w\in W}\sgn(w)e^{w(\rho+e(\g))})
$$
So both sides are stable (up to sign) under the Weyl group and hence it is sufficient to compare the 
coefficients of the terms corresponding to dominant weights, i.e. we have to 
prove for $\Om:=\{(w,\delta)\mid w\in W,\delta\in\G(\g_\lam),w(\rho)+e(\delta)\in  X^\vee_+\}$:
$$
\sum_{(w,\delta)\in\Om}\sgn(w)e^{w(\rho)+e(\delta)}
=\sum_{\delta\in \G(\g_\lam,{\rm dom})}e^{\rho+e(\delta)}.
$$
Let $\Om_0$ be the set of pairs $(w,\pi)\in\Om$ such that $w$ is the identity and 
$\delta\in \G(\g_\lam,{\rm dom})$. Set $\Om':=\Om-\Om_0$, to prove the proposition we have to show:
$$
\sum_{(w,\delta)\in\Om'}\sgn(w)e^{w(\rho)+e(\delta)}=0.\eqno{(*)}
$$
We will define an involution $\varphi:\Om'\rightarrow \Om'$
such that $\varphi(w,\delta)=(w',\delta')$ has the property: $\sgn(w)=-\sgn(w')$ and 
$w(\rho)+e(\delta)=w'(\rho)+e(\delta')$. This implies obviously $(*)$ and 
hence the proposition.

{\it The construction of the involution}: suppose first $(w,\delta)\in\Om'$
is such that $w$ is not the identity. Denote by $w(\rho)+\delta$ the gallery
of source $w(\rho)$ and target $w(\rho)+e(\delta)$ obtained by 
shifting the gallery $\delta$ by the co-character $w(\rho)$.
Since $w(\rho)+e(\delta)\in X^\vee_+$, the shifted gallery has to meet 
at least once a proper face of the dominant Weyl chamber $\Lc_f$, i.e., at least 
one of the small faces of $w(\rho)+\delta$ is contained in a proper
face of $\Lc_f$. If $w$ is the identity, then the shifted gallery also has to meet a proper 
face $F$ of $\Lc_f$, the pair would otherwise be an element of $\Om_0$. 

For a proper face $F$ of $\Lc_f$ denote by $\Om'(F)$ the set of pairs $(w,\delta)\in\Om'$ having the following property:  there exists a small 
face $w(\rho)+\Gamma_j'$ of the shifted gallery $w(\rho)+\delta$ 
such that the open faces $w(\rho)+{\Gamma_j'}^o\subset F^o$ are
contained in each other, and for all $k>j$  the open part
of the small faces $w(\rho)+{\Gamma_k'}^o$ are contained 
in the interior of $\Lc_f$.

The set $\Om'$ is obviously the disjoint union of the $\Om'(F)$,
so it is sufficient to define the involution for such an $\Om'(F)$.
Let $\al$ be a simple root orthogonal to $F$. For
$(w,\delta)\in\Om'(F)$ set $n:=\langle w(\rho),\al\rangle$, note that $n\not=0$. 

If $n<0$, then Lemma~\ref{simplelemmaA} ({\it v}) implies that $f_\al^{\vert n\vert}(\delta)$
is defined. Of course we have $w(\rho)+e(\delta)=s_\al w(\rho)+f^{\vert n\vert}_\al e(\delta)$,
and, by construction, 
$$
\varphi(w,\delta):=(s_\al w,f^{\vert n\vert}_\al(\delta))\in\Om'(F).
$$
Similarly, if $n> 0$, one sees that $e_\al^{n}(\delta)$
is defined, we have $w(\rho)+e(\delta)=s_\al w(\rho)+e^{n}_\al e(\delta)$,
and, by construction, 
$$
\varphi(w,\delta):=(s_\al w,e^{n}_\al(\delta))\in\Om'(F).
$$
Lemma~\ref{simplelemmaA} implies that $\varphi$
is an involution, which finishes the proof.
\endpf

Obviously, the character of $V(\lam)$ occurs in the decomposition above
with multiplicity one. We want to show that this character comes from the subset
of LS--galleries. 
\begin{rem}\label{foper}\rm
Let $\de\in \G(\g_\lam)$ be a combinatorial gallery having as target $\nu=e(\de)$ and
denote $t_{-\nu}$ the translation by $-\nu$. The combinatorial gallery $t_{-\nu}(\de)$ 
starts in $-\nu$ and ends in $0$. Let $\de^*$ be the gallery having the same faces 
as $t_{-\nu}(\de)$ but in reverse order, so the source of $\de^*$ is the origin, and 
the target is $-\nu$. The type of $\de^*$ is the same as the type of $\g_\lam^*$, 
which is the same as the type of $w_0(\g_\lam^*)$.  The latter is again a minimal 
gallery, but now joining the origin with the coweight $\lam^*=-w_0(\lam)$, the 
highest weight of the representation dual to $V(\lam)$. 
On sees easily that the map $\de\mapsto\de^*$ induces a bijection 
$\G(\g_\lam)\rightarrow  \G(w_0(\g_\lam^*))$, and, if $e_\al(\de)$
respectively $f_\al(\de)$ are defined, then
$$
e_\al(\de)=(f_\al(\de^*))^*\quad {\rm and}\quad f_\al(\de)=(e_\al(\de^*))^*
$$
\end{rem}
\begin{lem}\label{keepsposifolding}
\begin{itemize}
\item[{i)}] If $\g\in\Gamma(\g_\lam)$ and $e_\al\g$ (respectively $\widetilde{e}_\al\g$) is defined,
then $\dim e_\al\g=\dim\g +1$ (respectively $\dim \widetilde{e}_\al\g=\dim\g+1$),
and, similarly, if $f_\al\g$ is defined, then $\dim f_\al\g=\dim\g -1$.
\item[{ii)}] If $\delta \in \Gamma^+(\g_\lam)$ is such that 
$\widetilde{e}_\al\g$ is defined, then 
$\widetilde{e}_\al\g\in \Gamma^+(\g_\lam)$.
\item[{iii)}]  If $\delta \in \Gamma^+(\g_\lam)$ is such that $\wte_\al \delta$ is not defined but
$e_\al\delta$ (respectively $f_\al\delta$), then $e_\al\delta$ (respectively $f_\al\delta$)
is again positively folded.
\end{itemize}
\end{lem}
\proof
Consider first the operator $e_\al$. Translations do not affect the relative position 
of a pair $\Gamma_\ell'\subset \Gamma_\ell$ and $\Lc_{-\infty}$ with respect 
to a wall $\ch\supset \Gamma'_\ell$. Moreover, an affine
reflection associated to a simple root $\al$ affects the relative position only if
$\G_\ell'\subset\ch_{\al,p}$ for some $p\in\bz$. So to check
the dimension formula, it suffices to compare
the pair $\Gamma'_{j}\subset\Gamma_j$ (where $j$ is chosen as in 
Definition~\ref{defealpha}) with the corresponding pair of large and 
small faces in the new gallery. Now, by the definition of the operator,
$\ch_{\al,m+1}$ is not load-bearing for $\g$ at $\G_j$, but 
$\ch_{\al,m+1}$ is load-bearing for $e_\al\g$. The same arguments
apply to the operator $f_\al$, only in this case $\ch_{\al,m}$ is load-bearing for
$\g$ at $\G_j$ but not  for $f_\al\g$. 

Using the same arguments as above, one sees that for the operator $\wte_\al$
it suffices to consider the pairs $\Gamma'_{j}\subset\Gamma_j$ and 
$\Gamma'_{k}\subset\Gamma_k$ ($j,k$ chosen as in Definition~\ref{defetildealpha}).
Now $\ch_{\al,m}$ is load-bearing for $\g$ and $\wte_\al \g$ at $\G_k$, and
$\ch_{\al,m}$ is load-bearing for $\wte_\al\g$ at $\Gamma_j$ but not for $\g$.
This implies the dimension formula. 

To prove that $\wte_\al\delta \in \Gamma^+(\g_\lam)$ if 
$\delta \in \Gamma^+(\g_\lam)$, arguments similar to the ones 
above reduce the proof to check what the operator changes 
for the triples  $\Gamma_{j-1}\supset\Gamma'_{j}\subset\Gamma_j$ 
and $\Gamma_{k-1}\supset\Gamma'_{k}\subset\Gamma_k$. 
But since the additional foldings induced by $\wte_\al$ at $\Gamma'_{j}$ and 
$\Gamma'_{k}$ are positive by construction, it follows that 
$\widetilde{e}_\al\g\in \Gamma^+(\g_\lam)$.

It remains to prove {\it iii)}. Note that $\wte_\al \delta$ is not defined
if and only if $\wte_\al \delta^*$ is not defined. Since $\delta$ is positively folded if and only if 
$\delta^*$ is, by Remark~\ref{foper} 
it suffices to give a proof for the operator $e_\al$.
Once again, arguments as above show that it suffices to check what the
operator changes for the triples  $\Gamma_{j-1}\supset\Gamma'_{j}\subset\Gamma_j$ 
and $\Gamma_{k-1}\supset\Gamma'_{k}\subset\Gamma_k$
($j,k$ chosen as in Definition~\ref{defealpha}). The operator
$e_\al$ adds an additional positive folding at $\Gamma_{j}'$, which finishes the
proof for the index $j$.

The assumption that $\wte_\al\delta$ is not defined implies that 
$\Gamma_{k-1}$ and $\Gamma_k$ are on the same side of $\ch_{\al,m}$. Let us denote by $\De_{k-1} \supset \De'_k\subset \De_k$ the triple in $e_\al\de$ corresponding to $\G_{k-1}\supset \G'_k\subset \G_k$. If $\lam$ is regular then the two alcoves $\De_{k-1}$ and $\De_k$ are separated by the wall spanned by $\De'_k$, namely $\ch_{\al, m+2}$. If $\lam$ is not regular, the two large faces $\De_{k-1}$ and $\De_k$ are still separated by $\ch_{\al, m+2}$. And, because $\al$ is simple, the possible other positive foldings at the place $k$ in $\de$ are sent by $e_\al$ to positive foldings.

It follows that $e_\al\delta\in\G^+(\g_\lam)$.
\endpf

We recall now Proposition~\ref{dim-inequality} from section~\ref{dimensiongallery}
and give the proof:

\par\vskip 5pt\noindent
{\bf Proposition~\ref{dim-inequality}.}
{\it If $\delta\in\Gamma^+(\gamma_\lam,\nu)$, then $\dim\de\le \langle\lam+\nu,\rho\rangle$.}
\par\vskip 5pt
\proof
We proceed by induction on ${\rm ht}(\lam-\nu)$, the height of $\lam-\nu$. Note that 
$\langle\lam+\nu,\rho\rangle=2\langle\lam,\rho\rangle-{\rm ht}(\lam-\nu)$.
There is only one gallery $\delta\in \G^+(\g_\lam)$ such that $e(\delta)=\lam$,
it is $\delta=\gamma_\lam$. The dimension formula holds in this case by 
Example~\ref{exampleextremalgallery}.

If $\nu=e(\delta)\not=\lam$, then $\delta=[\de_0,\ldots,\de_p]$ for some $\de_0 \not= id$.
Let $\al$ be a simple root such that $s_\al \de_0 < \de_0$, then
either the gallery $\delta'=\wte_\al\delta$ is defined and is positively folded,
or $\wte_\al\delta$ is not defined, but then the gallery $e_\al\delta$ is defined
and positively folded (Lemma~\ref{keepsposifolding}). In the latter case the
claim follows immediately by induction and Lemma~\ref{keepsposifolding} {\it i)}.
In the first case we repeat the procedure (if necessary, with another simple root)
with $\delta'$ until we get a positively folded gallery $\delta''$ such that still $\nu=e(\delta'')$,
$\wte_{\al'}\delta'$ is not defined but $e_{\al'}\delta''$ is. Again, the claim
follows now by induction and Lemma~\ref{keepsposifolding} {\it i)}.
\endpf

As a consequence of the proof above we get, 
in particular, a proof of Proposition \ref{characterformula}:
\begin{coro}\label{lsgalleryandeal}
\begin{itemize}
\item[({\it i\/})]  $\dim\delta=\la\lam+\mu,\rho\ra$ 
if and only if there exists a sequence of simple roots
such that $e_{\al_{i_1}}\cdots e_{\al_{i_t}}\delta=\gamma_\lam$, or, equivalently,
$\delta=f_{\al_{i_t}}\cdots f_{\al_{i_1}}\gamma_\lam$. In particular, the set of
LS--galleries in $\Gamma^+(\g_\lam)$ is the subset
generated from $\gamma_\lam$ using the
operators $f_\al$.
\item[({\it ii\/})] The gallery $\de$ is an LS--gallery if and only if $\de^*$ is an LS--gallery.
\item[({\it iii\/})] $\chara V(\lam)=\sum \exp(e(\g))$, where the sum runs over all
LS--galleries in $\G_{LS}^+(\g_\lam)$.
\end{itemize}
\end{coro}
\proof
If $\delta$ is an LS--gallery, then the proof above shows
(for dimension reasons) the existence of a sequence such that 
$e_{\al_{i_1}}\cdots e_{\al_{i_t}}\delta=\gamma_\lam$.
Conversely, if such a sequence exists, then we get equivalently a 
sequence such that $\delta=f_{\al_{i_t}}\cdots f_{\al_{i_1}}\gamma_\lam$.
Suppose we  have already shown that 
$\delta'=f_{\al_{i_j}}\cdots f_{\al_{i_1}}\gamma_\lam$, $0\le j<t$,
is an LS--gallery. Then Lemma~\ref{keepsposifolding} {\it i)} and Proposition~\ref{dim-inequality}
implies that $\wte_\al\delta'$ is not defined (for dimension reasons), so $ f_{\al_{i_{j+1}}}\delta'$
is positively folded and an LS--gallery by Lemma~\ref{keepsposifolding} {\it i)} and {\it iii)}.
It follows by induction that $\delta$ is an LS--gallery. Since the set of LS--galleries
is closed under the folding operators, the proof of Corollary~\ref{characterdecompo}
applies also to the set of LS--galleries in $\G^+(\g_\lam)$. Since $\g_\lam$ is the
only dominant LS--gallery, this proves part {\it iii)}. Part {\it ii)} follows from 
Remark~\ref{foper}, part {\it i)} and the fact that $\gamma_\lam^*$ is an LS--gallery.
\endpf

We will now describe the connection between the path model and the LS--galleries.
Recall that the minimal galleries are of the form $[\g_0,\tau^m_1,...,\tau^m_p]$, 
where $\g_0\in W$ is a representative of minimal length of a class in $W/W_\lam$, 
and each $\tau^m_j\in W'_j$  is the unique minimal representative of the class of 
the longest element in $W'_j$.
\begin{dfn}\rm
The {\it companion $(\sigma_0,\ldots,\sigma_p)$ of a gallery} $\g\in\G(\g_\lam)$ 
is a sequence of Weyl group cosets in $W/W_\lam$ satisfying the following 
conditions:

$\bullet$ $\sigma_0 =\gamma_0 \mod W/W_\lam$ 

$\bullet$ if there is no folding at $\Gamma'_j$, then $\sigma_j=\sigma_{j-1}$ 

$\bullet$  assume that  $\g$ is folded at $\Gamma'_j$;
\par\noindent
{\it a)} If $\lam$ is regular, then let $s_{\beta}$ be the reflection used for the folding, 
i.e., $\Gamma_j=s_{\beta,m}\Gamma_{j-1}$ for an appropriate integer $m$. 
We set $\sigma_{j}=s_{\beta}\sigma_{j-1}$ in $W$.
\par\noindent
{\it b)} If $\lam$ is not regular, then more than one reflection can occur in a folding. Let
$s_{\beta_1},\ldots,s_{\beta_r}$ be the reflections used for the folding, i.e.,
$\Gamma_j=s_{\beta_1,m_1}\cdots s_{\beta_r,m_r}\Gamma_{j-1}$ for appropriate
integers $m_1,\ldots,m_r$, then set 
$\sigma_{j}=s_{\beta_1}\cdots s_{\beta_r}\sigma_{j-1}$ in $W/W_\lam$.
\end{dfn}

\begin{rem}\label{companionandpathrem}\rm
To compare this notion with the language of the path model, note that one
can find a piecewise linear path $\pi$ starting in $0$ and ending in $\lam$,
such that all directions are dominant rational weights and have the same 
stabilizer in $W$ as $\lam$, and the image is contained in 
$\g_\lam=(F_f\subset \G_0\supset\G_1'\subset\G_1\supset\ldots\subset\G_p\supset F_\lam)$. 
\end{rem}
\begin{exam}\label{companionandpathexam}\rm
Assume for simplicity that $G$ is simple. Take the piecewise linear path that joins 
the origin with the barycenter of $\G_1'$, then the barycenter of $\G_1'$ with the 
barycenter of $\G_2'$, etc, and finally joins the barycenter of $\G_p'$ with $\lam$. 
Note that the segments in $\Gamma_i$, $1\le i\le p-1$, are parallel to all codimension 
one faces of the $\G_i$ except for the faces $\G_i'$ and $\G_{i+1}'$. So the minimality 
of the gallery implies that each segment is a positive rational multiple of a dominant 
weight. Now these may have a larger stabilizer than $\lam$, but at least the first segment 
has obviously the desired property. Now by moving around the turning point on $\G_1'$ 
(and letting the others fixed), we may still assume that the first segment is dominant rational
and has the same stabilizer as $\lam$,  and we may assume in addition that this is now 
also true for the second segment. Continuing the finite procedure, we obtain the desired 
piecewise linear path.

Now the image $\pi_\g$ of this path in a folded gallery $\g$ of the same type will be a path 
such that its directions will be Weyl group conjugates of the directions we
started with. The corresponding Weyl group elements (or rather the classes) 
give the companion.
\end{exam}

\begin{dfn}\rm
A gallery is called a {\it combinatorial LS--gallery} if 
the following holds: for each pair $(\sigma_{i-1},\sigma_i)$ in the companion, 
the sequence of folding reflections $s_{\beta_1,m_1},\ldots,s_{\beta_r,m_r}$
can be chosen such that applied to $\sigma_{i-1}$, this gives a 
sequence $\sigma_{i-1} =\tau'_0 > \tau'_1=s_{\beta_1}\tau'_0 >\ldots>\tau'_r=
s_{\beta_r}\tau_{r-1}=\sigma_i$  
of Weyl group cosets (modulo $W_\lam$) such that  the length is always decreasing 
by one for each reflection.
\end{dfn}

\begin{rem}\rm
The condition that the classes are decreasing in the Bruhat order implies that
this is a positive folding: if $\kappa=s_\beta\tau>\tau$, then $\kappa(\lam)\prec\tau(\lam)$
in the weight ordering (i.e., $\tau(\lam)-\kappa(\lam)$ is a non-negative sum of positive roots),
and hence $\la\kappa(\lam),\beta\ra>0$.
\end{rem}

The special r\^ole played by the combinatorial LS--galleries $\g$ is the following: let $\al$
be a simple root and let $m\in\bz$ be minimal such that a small face of $\g$ is contained 
in $\ch_{\al, m}$.
\begin{lem}
If $\g$ is a combinatorial LS--gallery and $e_\al\g$ (respectively $f_\al\g$) is well defined, then 
$e_\al\g$ (respectively $f_\al\g$) is again a combinatorial LS--gallery of type $t_{\g_\lam}$. 
\end{lem}
\proof
Suppose first that $\wte_\al\g$ is not defined.
We know already (Lemma~\ref{keepsposifolding} {\it iii)}) 
that if, in this case, $e_\al\g$ is defined, then it is a positively folded gallery.
To see that it is in fact a combinatorial LS--gallery, one can now apply the same arguments
as in \cite{Lit1} on chains of Weyl group cosets to prove that the companion of $e_\al\g$
has again the desired properties. The proof for $f_\al$ follows in the same
way using the ``$*$''--operation.

It remains to show that $\wte_\al\g$ is not defined for a combinatorial
LS--gallery, i.e., we have to show that the gallery does not cross the hyperplane
$\ch_{\al,m}$. Let $\pi$ be a piecewise linear path corresponding to $\gamma_\lam$
as in Remark~\ref{companionandpathrem} and Example~\ref{companionandpathexam},
and denote $\pi_\g$ the corresponding folded path for $\g$. We may assume that the
path has rational turning points and is linear inside each large face of the gallery. 
Suppose that $\wte_\al\g$ is defined, i.e., the gallery crosses $\ch_{\al,m}$,
let $\G'_j\subset \G_{j}$ and $\G_{k-1}\subset \G_{k}'$ be the small and large faces
as in Definition \ref{defetildealpha}. Let $\mu_j,\mu_{k-1}$ be the rational coweights
that correspond to the segments of $\pi_\g$ in the large faces $ \G_{j},  \G_{k-1}$.
Since $\g$ crosses $\ch_{\al,m}$ we have $\langle\mu_j,\al\rangle<0$, and
since $k>j$ is minimal such that $\G_{k}'\subset\ch_{\al,m}$, we have 
$\langle\mu_{k-1},\al\rangle>0$. Recall that none of the folding hyperplanes
at the small faces $\G'_{j+1},\ldots, \G'_{k-1}$ is of the form $\ch_{\al,p}$.
But then such a change of sign is impossible for a combinatorial LS--gallery by the
following Lemma~\ref{rootsignswitchlem} and Remark~\ref{rootsignswitchrem}, 
which finishes the proof.
\endpf

Let $\mu,\nu$ be dominant rational coweights having the same stabilizer $W'$
in $W$. For $\sigma\in W/W'$ denote $N(\sigma)$ the set of positive roots $\beta'$ such that 
$\langle\sigma(\mu),\beta'\rangle\le 0$. Let $\beta\in N(\sigma)$ be such that
$\ell(s_\beta\sigma)=\ell(\sigma)-1$ in $W/W'$. 
\begin{lem}\label{rootsignswitchlem}
If $\al$ is a simple root
such that $\langle\sigma(\mu),\al\rangle<0$ but $\langle s_\beta\sigma(\nu),\al\rangle>0$,
then $\al=\beta$.
\end{lem}
\begin{rem}\label{rootsignswitchrem}\rm
Note that the positive roots
occurring in the folding of a combinatorial LS--gallery $\g$ have exactly this property:
let $\Delta_{l-1}\supset\Delta_l'\subset\Delta_l$ be two large and one small face
such that $\Delta_l'\subset \ch_{\beta,n}$, and the affine hyperplane separates
the large faces. Suppose that all foldings at the
small faces $\Delta'_r$, $r\le l$ satisfy the condition for a combinatorial LS--gallery.
Let $\sigma_l$ be the element in the companion for $\Delta'_l$, and let $\mu_l$ 
be the segment in $\De_{l-1}$ of an associated path $\pi_\g$ as in the proof above. 
Because any combinatorial LS--gallery is positively folded, we need $\langle \mu_l,\beta\rangle<0$, and to satisfy the condition for a combinatorial
LS--gallery, we also need $\ell(s_\beta\sigma_l)=\ell(\sigma_l)-1$.
\end{rem}
\proof
The following proof is a slightly adapted version of the proof of 
Proposition~1 in \cite{Dab}, so we will skip some details.  
For a given $\beta\in N(\sigma)$, instead of decomposing
$N(\sigma)$ into two sets, one decomposes it into three sets. Let
$E(\sigma):=\{\theta\in N(\sigma)\mid \langle\sigma(\mu),\theta\rangle
=0, s_\beta\theta\prec 0\}$ and
$$
A(\sigma):=\{\theta\in N(\sigma)\mid s_\beta\theta\succ 0\},\ 
B(\sigma):=\{\theta\in N(\sigma)\mid \langle\sigma(\mu),\theta\rangle
<0, s_\beta\theta\prec 0\},
$$
then $N(\sigma)$ is the disjoint union of $A(\sigma)$, $B(\sigma)$ and $E(\sigma)$.
One checks easily that $A(s_\beta\sigma)=s_\beta A(\sigma)$, 
$E(s_\beta\sigma)=-s_\beta E(\sigma)$ and $B(s_\beta\sigma)\cup\{\beta\}\subset
B(\sigma)$. One shows then: if $\ell(s_\beta\sigma)=\ell(\sigma)-1$, then 
$E(\sigma)=\emptyset$ and $B(\sigma)=B(s_\beta\sigma)\cup\{\beta\}$.

Let $\al,\beta$ be as above and suppose $\al\not=\beta$. If 
$s_\beta(\al)\succ 0$, then $s_\beta(\al)\in A(s_\beta\sigma)$. 
But note that $\al,\beta\in N(\sigma)$ implies $\al+k\beta\in C(\sigma)$, the real convex cone spanned by $N(\sigma)$. Since $s_\beta(\al)\succ 0$, this shows $s_\beta(\al)\in N(\sigma)$
and hence $s_\beta(s_\beta(\al))=\al\in N(s_\beta\sigma)$, in contradiction 
to the assumption.
So $s_\beta(\al)\prec 0$, and hence $\al\in B(\sigma)-\{\beta\}=B(s_\beta\sigma)$,
in contradiction to the assumption.
\endpf

Using the operators $e_\al$, we can transform a LS--gallery $\g$ into a LS--gallery $\g'$
having $\sigma_0=1$ for the companion, and hence  $\g'=\g_\lam$. It follows:
\begin{coro}
The set of LS-galleries in $\Gamma^+(\g_\lam)$ coincides with the set of combinatorial LS-galleries.
\end{coro}
For the gallery $\gamma_\lam$ let $\pi:[0,1]\rightarrow \ca$ be a piecewise linear path 
as in Remark~\ref{companionandpathrem}.
By comparing the definition of the folding operators $e_\al,f_\al$ on galleries
with the  the definition of the root operators $e_\al,f_\al$ on piecewise 
linear paths, one sees easily that the set of LS--galleries in $\G^+(\gamma_\lam)$ 
is exactly the set of
galleries obtained in the following way: let $B(\pi)$ be the path model for the representation
$V(\lam)$ obtained by applying the root operators to $\pi$. To a path $\eta\in B(\pi)$ we
associate the gallery $\gamma_\eta$ defined as the sequence of faces 
gone through by $\eta$. Then $\G^+_{LS}(\gamma_\lam)=\{\gamma_\eta\mid\eta\in B(\pi)\}$,
and this identification is equivariant with respect to the operators $e_\al,f_\al$.
As a consequence we get (by \cite{Lit2}, \cite{J}, \cite{kashsim}):
\begin{thm}\label{crystal}
Let $B(\gamma_\lam)$ be the directed colored graph having as vertices the set of LS-galleries in 
$\Gamma^+(\gamma_\lam)$, and put an arrow $\delta\mapright{\alpha}\delta'$
with color $\alpha$ between two galleries if $f_\alpha(\delta)=\delta'$. Then this graph is
connected, and it is isomorphic to the crystal graph of the irreducible 
representation $V(\lam)$ of $G^\vee$ of highest weight $\lam$.
\end{thm}


\section{Variety of Galleries}\label{varietiesofgalleries}

From now on, we take $G$ to be simply connected. By Remark \ref{orbits}, this assumption is not a restriction. In this section, we give two equivalent definitions of the Bott-Samelson variety, and 
we point out the first step of the connection between galleries and MV-cycles.
Let $\lam\in \Xveep$ be a dominant co-character and let $t_\lam = type (F_\lam)$ be 
the type of the vertex of $\tta$ associated to $\lam$. Let $\cq_\lam$ be the parahoric 
subgroup of type $t_\lam$ containing $\cb$. The Bott-Samelson variety $\BSlam$
occurs as a natural resolution (see \cite{Ku}) of the "Schubert variety" $X_\lam = \overline{\cg_\lam}$, where $\cg_\lam =
G(\co)\lam\hookrightarrow G(\ck)/\cq_\lam$ (or $\cg_\lam = \cp_\co\lam\hookrightarrow\hat\cl(G)/\hat\cq_\lam$ in terms of the Kac-Moody group, see section \ref{KMgroup}). By a resolution we mean that $\BSlam$ is a smooth algebraic 
variety and the morphism $\pi : \BSlam\to X_\lam$ is birational and proper.

Let 
$$
t_{\g_\lam} = (S = t'_0\supset t_0 \subset t'_1 \supset\cdots\supset
t_{j-1} \subset t'_j \supset t_j \subset\cdots \subset t'_p \supset t_p \subset t_{\lam}).
$$ be the gallery of types of the fixed minimal combinatorial gallery $\g_\lam$ 
joining $F_f$ with $F_\lam$. Consider 
first the usual definition of $\BSlam$ in terms of $\hat\cl(G)$ (in \cite{Ku} $\S$7.1, S. Kumar gives the 
definition only for the case where 
the $t_j$'s are trivial, but it makes also sense in our case). For $0\leq j\leq p$, let us 
denote by $\hat\cp_j$ (resp. $\hat\cq_j$) the parabolic subgroup 
of type $t'_j$ (resp. $t_j$) containing $\hat\cb$, of course $\cp_0 = \cp_\co$.  

\begin{dfn}\label{defiBott1}
The Bott-Samelson variety $\BSlam$ is defined as
$$
\BSlam = \cp_\co\times_{\hat\cq_0}\hat\cp_1\times_{\hat\cq_1}\cdots\times_{\hat\cq_{p-1}}\hat\cp_p/{\hat\cq_p},
$$ 
i.e. the algebraic (complex) variety defined as the quotient of the group 
$\cp_\co\times \hat\cp_1 \times\cdots\times \hat\cp_p$ by the subgroup 
$\hat\cq_0\times \hat\cq_1\times\cdots\times \hat\cq_p$ under the (right) action given by 
$g \cdot q = (g_0q_0,q_0^{-1}g_1q_1,...,q_{p-1}^{-1}g_pq_p)$ 
where $q = (q_0,q_1,...,q_r) \in \hat\cq_0\times \hat\cq_1\times\cdots\times \hat\cq_p$ 
and $g = (g_0,g_1,...,g_{p})\in \cp_\co\times \hat\cp_1 \times\cdots\times \hat\cp_{p}$. 
\end{dfn}

The Bott-Samelson variety $\BSlam$ is a 
smooth, projective (complex) variety of dimension $\dim\ (\BSlam) = 2\la\lam,\rho\ra$.  
We will denote by $g = [g_0,g_1,...,g_p]$ the points of this variety and we call them {\it galleries}
(the reason will become apparent soon). The Bott-Samelson variety naturally comes also equipped 
with a proper and birational morphism 
$$
\pi : \BSlam\to X_\lam
$$ 
by associating to a gallery $g = [g_0,g_1,...,g_p]$ the class of the product $g_0g_1\cdots g_p$ 
in the affine Grassmannian $\hat\cl(G)/\hat\cq_\lam$. 
The following Definition--Proposition is a slight generalization of a result of C. Contou-Carr\`ere $\cite{CC}$.

\begin{defiprop}\label{defipropbott}
For $0\leq j\leq p$, let us 
denote by $\cp_j$ (resp. $\cq_j$) the parahoric subgroup of $G(\ck)$ 
of type $t'_j$ (resp. $t_j$) containing $\cb$.  
The Bott-Samelson variety $\BSlam$ is the closed subvariety of the product
$$
G(\ck)/\cq_0\times\Big (\prod_{1\leq j\leq p} G(\ck)/\cp_j\times G(\ck)/\cq_j\Big )\times G(\ck)/\cq_\lam
$$
given by all the sequences of parahoric subgroups of the shape 
$$
(G(\co)\supset \cq'_0\subset \cp'_1 \supset \cq'_1 \subset \cdots \subset \cp'_p\supset \cq'_p\subset \cq'_\lam),
$$
where $type(\cp'_i) = t'_i$, $type(\cq'_j) = t_j$ and $\cq'_\lam$ is a subgroup 
associated to a vertex of type $t_\lam$. In other words, the set of points of $\BSlam$ 
is given by the set of all the galleries in the affine building $\cj^{\mathfrak a}$ of type 
$t_{\g_\lam}$ starting at $F_f$ : 
$$
g = (F_f\subset F_0\supset F'_1 \subset F_1 \supset \cdots \supset F'_p\subset F_p\supset F'_\lam),
$$
where the type of each face corresponds to the type of the corresponding parahoric subgroup. 
\end{defiprop}

\proof
In \cite{CC}, Contou-Carr\`ere proves the proposition in the case of a Bott-Samelson variety built of parabolic 
subgroups of a reductive group scheme, but it applies also to the Kac-Moody case. Let us recall that the Kac-Moody group acts on the building. The equivalence between the two definitions is then obtained by the 
isomorphism which associates to a point $[g_0,g_1,...,g_p]$  in $\BSlam$ the gallery 
$$
g = (F_f\subset F_0\supset F'_1 \subset F_1 \supset \cdots \supset F'_p\subset F_p\supset F'_{p+1})
$$ 
where, for $0\leq j\leq p$, $F_j = g_0g_1\cdots g_jF_{\cq_j}$ ($F_{\cq_j}$ is the only face of 
the fundamental alcove $\De_f$ of type $t_j$). The faces $F_j'$ of the gallery $g$ are uniquely 
determined by the $F_j$'s and the type $t_{\g_\lam}$. 
\endpf

\begin{rem}\label{gallerypoints}\rm
The identification of galleries of type $t_{\g_\lam}$ and of source $F_f$ with the points in the Bott-Samelson
variety allows us to view the combinatorial galleries $\G(\g_\lam)$ of type $t_{\g_\lam}$ and of source $F_f$
as points in $\BSlam$, they are precisely the $T$--fixed points.
Here the action of $G(\co)$ on the galleries is the action induced
by the operation on the building. Since $G(\co)$ operates on $\cj^{\mathfrak a}$ by simplicial maps, 
the action preserves the type of a gallery, and, if the gallery is minimal, then so are all galleries in the orbit. 
\end{rem}

Further, in this setting, the morphism $\pi : \BSlam\to X_\lam$ $\pi$ is just the restriction of the projection on the last factor. Hence, it maps the minimal gallery $\gamma_\lam$ 
onto the  point $\lam$, so $\pi$ induces a morphism between $G(\co)\gamma_\lam$ and
$\cg_\lam$. Now one sees easily that all minimal galleries of type  $t_{\g_\lam}$ are in one
$G(\co)$--orbit, so by the birationality of $\pi$ we get:

\begin{lem}\label{openorbitgallery}
The desingularization map $\pi : \BSlam\to X_\lam$ induces a bijection between 
the set of all minimal galleries $\hat\cg_\lam$ of type $t_{\g_\lam}$ in $\BSlam$ and the open orbit 
$\cg_\lam\subset X_\lam$.
\end{lem}

The retraction  $r_{-\infty}$ extends naturally to galleries by applying  $r_{-\infty}$ simultaneously
 to all faces in the gallery. Since  $r_{-\infty}$ is a map of simplicial complexes, it preserves the type,
 i.e., the image of a gallery of type $t_{\gamma_\lam}$ is a combinatorial gallery of type $t_{\gamma_\lam}$. 
Let $\G(\g_\lam)$ be the set of all combinatorial galleries of type $t_{\gamma_\lam}$, i.e., the
set of $T$-fixed points (Remark~\ref{gallerypoints}) in $\BSlam$.
The retraction map on the galleries can be described geometrically as follows:

\begin{prop} \label{BBcell}
The retraction with center at $-\infty$ induces a map
$\hat{r}_{\gamma_\lambda}:\BSlam\to\Gamma(\gamma_{\lam})$.
The fibres $C(\delta)=\hat{r}_{\gamma_\lambda}^{-1}(\de)$, $\de\in\G(\g_\lam)$, are naturally endowed
with the structure of a locally closed subvariety, each of them being isomorphic to an affine space.
In fact, the $C(\delta)$ are the Bialynicki-Birula cells $\{x\in \BSlam\mid \lim_{s\rightarrow 0}  \eta(s)x=\delta  \}$

of center $\delta$ in $\BSlam$, associated to a generic one--parameter
subgroup $\eta$  of $T$ in the anti--dominant Weyl chamber.
\end{prop}

\proof
The map $\hat{r}_{\gamma_\lambda}$ is well--defined, it remains to describe the fibres.
Let $\de$ be a combinatorial gallery in $\G(\g_\lam)$, $\de= (F_f=\Delta_0'\subset \Delta_0\supset\Delta'_1
\subset\cdots\subset\De_p \supset \De_{p+1}=F_\nu)$, and suppose that
$g = (F_f=F_0'\subset F_0\supset  \cdots \supset F'_p\subset F_p\supset F'_\lam)$
is a gallery of type $t_{\gamma_\lam}$ which retracts onto $\de$. By Proposition~\ref{fibresofrinf},
we know that there exist hence $u_i\in U^-(\ck)$, $i=0,\ldots,p$, such that
$F_i=u_i\De_i$, and, since $g$ is a gallery (of type $t_{\gamma_\lam}$),
we have necessarily  $u_0\De'_0=\De'_0$ and $u_{i-1}\De_i'=u_i\De_i'$,
$i=1,\ldots,p$. Conversely, given any sequence $(u_0,\ldots,u_p)$ of elements
in $U^-(\ck)$ satisfying $u_{i-1}^{-1} u_i\in \hbox{\rm Stab}_{G(\ck)}(\Delta_i')$, the stabilizer
of the face in $G(\ck)$,
for $i=0,\ldots,p$ (where we set $u_{-1}=1$), then
\begin{equation}\label{ugallery}
g = (F_f=\Delta_0'\subset u_0\Delta_0\supset u_1\Delta'_1\subset u_1\De_1\supset
\cdots\supset u_p\De'_p\subset u_p\De_p \supset u_p \De_{p+1}=u_p F_\nu)
\end{equation}
is a gallery of type $t_{\gamma_\lam}$, which retracts onto $\de$. The maximal torus $T$ normalizes
$U^-(\ck)$ and  $(tu_{i-1}t^{-1})(tu_i^{-1}t^{-1})\in \hbox{\rm Stab}_{G(\ck)}(\Delta_i')$
for $t\in T$ if and only if this holds for $u_{i-1} u_i^{-1}$, so the fibres of the retraction
map $\hat{r}_{\gamma_\lambda}$ are stable under the natural $T$--action on $\BSlam$.
Let $\eta:\bc^*\rightarrow T$ be a one--parameter subgroup in generic position in
the open anti-dominant Weyl chamber, so $\lim_{s\rightarrow 0}(\eta(s)u\eta(s)^{-1})=1$ for $u\in U^-(\ck)$.
If $g$ is as in (\ref{ugallery}), then $\lim_{s\rightarrow 0}\eta(s)g=\delta$ because the
$\De_i,\De_i'$ are $T$-fixed points in the building, and hence
$$
\begin{array}{rcl}
\eta(s)g&=&(\eta(s)F_f=\eta(s)\Delta_0'\subset \eta(s)u_0\Delta_0\supset \eta(s)u_1\Delta'_1\subset
\cdots\supset \eta(s) u_p\De_p \supset \eta(s)u_p F_\nu)\\
&=&(F_f\subset (\eta(s)u_0\eta(s)^{-1})\Delta_0\supset (\eta(s)u_1 \eta(s)^{-1})\Delta'_1\subset \cdots\\
&&\hfill\cdots\supset (\eta(s) u_p\eta(s)^{-1})\De_p \supset (\eta(s)u_p \eta(s)^{-1}) F_\nu),\\
\end{array}
$$
which implies the claim. The $T$--fixed points in $\BSlam$ are the combinatorial galleries
in $\G(\g_\lam)$. Since $\eta$ is generic, we can assume that the
$T$--fixed points and the $\eta$--fixed points coincide. So, by \cite{BB1} and \cite{BB2},
the subsets
$$
C(\delta)=\{x\in \BSlam\mid \lim_{s\rightarrow 0}  \eta(s)x=\delta  \}
$$
form a decomposition of $\BSlam$ into locally closed subsets, each of which is isomorphic
to an affine space, and the calculation above shows that the $C(\de)$ are the fibres of $\hat{r}_{\gamma_\lambda}$.
\endpf


By Lemma~\ref{openorbitgallery}, the object of our main interest, the orbit $\cg_\lam$,
can be identified  with the minimal galleries $\hat\cg_\lam$ in $\BSlam$. Here the special r\^ole of the 
positively folded galleries and its connection with the MV-cycles becomes apparent:

\begin{thm}\label{firstpartmaintheorem}
The restriction of $\hat{r}_{\gamma_\lam}$ induces a map 
$r_{\gamma_\lam}:\cg_\lam\to\Gamma^+(\gamma_{\lam})$.
Further, the union $\bigcup_{\delta\in\Gamma^+(\gamma_{\lam},\nu)}
r_{\gamma_\lam}^{-1}(\delta)$ of the fibres over the galleries in
$\Gamma^+(\gamma_{\lam})$ with target $\nu$
is the (set theoretic) intersection $U^-(\ck)\nu\cap\cg_\lam$.
\end{thm}

\proof
By Proposition~\ref{fibresofrinf}, we know that the fibres of 
$r_{-\infty}:\cj^{\mathfrak a}\rightarrow \ca$ are the $U^{-}(\ck)$--orbits.
If $\delta\in \BSlam$ is a point (or rather a gallery), then $r_{-\infty}(\pi(\delta))$
is the target of the gallery $\hat{r}_{\gamma_\lam}(\delta)$. So the union
of the fibres $r_{\gamma_\lam}^{-1}(\delta)$ with $\delta\in \Gamma(\gamma_{\lam},\nu)$
is precisely the set theoretic intersection of the orbit $U^-(\ck)\nu$ with
the orbit $\cg_\lam$.

To finish the proof of the theorem, it remains to show that is suffices
to consider positively folded galleries, but this follows from the next lemma.
\endpf
\begin{lem}
If $\de = [\de_0,\de_1,...,\de_p] = (F_f\subset \Delta_0\supset\Delta'_1\subset
\cdots\Delta_{j-1}\supset\Delta'_{j}\subset\Delta_{j}\supset\cdots  \supset F_\nu)$ 
is not a positively folded gallery, then $\hat{r}^{-1}_{\g_\lam}(\de)\cap\pi^{-1}(\cg_\lam) = \emptyset$. 
\end{lem}

\proof
The assumption means that there exists an index $j\in\set{1,...,p}$ and a 
wall $\ch\in\sh_j$ such that $\de_j \not = \tau^m_j$ and the image of the 
folding around $\ch$ is not separated from $\Lcinf$ by $\ch$. Particularly, 
the large face $\Delta_{j}$ is not separated from $\Lcinf$ by $\ch$ and, 
because we fold around $\ch$, the two consecutive large faces $\De_{j-1}$ 
and $\De_j$ are on the same side of this wall. Suppose that $j$ is the first of its kind.

Let $g\in\hat{r}^{-1}_{\g_\lam}(\de)$
be a gallery that retracts on $\de$. We can assume that $g$ has already been 
retracted onto $\de$ up to the index $j-1$, meaning that the 
gallery is of the form
$$
g =  [\de_0, \de_1,...,\de_{j-1},g_{j},...,g_p] = 
(F_f\subset \Delta_0\supset\cdots\subset\Delta_{j-1}\supset \Delta'_{j} 
\subset\ F_j\supset F'_{j+1} \cdots  \supset F_{p+1}),
$$
with $F_{j} = \de_0\de_1\cdots\de_{j-1} g_{j}F_{\cq_j}$, where $g_j\in\hat\cp_j/\hat\cq_j$ 
and $F_{\cq_j}$ is the face of type $t_j$ contained in the fundamental alcove. 

Now let us assume that $g$ is minimal, then $g^j = (\De_{j-1}\supset \De_j' \subset F_j)$ is also minimal. 
That implies that $\De_{j-1}$ and $F_j$ are separated by the wall image of $\ch$ in any apartment $A$ 
containing $g^j$. But, incarnating the retraction $\rinf$ using a suitable far away alcove $\Sigma\in A\cap\tta$ 
on the same side of $\ch$ than $\De_{j-1}$, we know that the retraction $\rinf = r_{\Sigma,\tta}$ 
preserves the distances from $\Sigma$ and reduces to an isomorphism of complex chambers 
$A\simeq\tta$. So, $\De_{j-1}$ and $\rinf (F_j)$ are also  separated by the wall $\ch$, which is 
not possible if we want $\rinf (F_j) = \De_j$. Hence, the lemma is proved.
\endpf


\section{An open covering of $\BSlam$}\label{opencovering}
The open subsets of the covering of $\BSlam$ will be centered at the combinatorial galleries $\G(\g_\lam)$.
This covering is a slight generalization of the one considered by J. Tits \cite{T82}. The open sets will be 
built using the root subgroups of the Kac-Moody group $\hat\cl(G)$. 

\begin{dfn}\rm
For any {\it real root} $\eta$ in the root system of the Kac-Moody group 
$\hat{\mathcal L}(G)$, there exists a one-dimensional root subgroup $\cu_\eta$,
isomorphic to $(\bc, +)$ (see \cite{Ku}, 6.2.7) having as Lie algebra the corresponding
root subspace in $\hat{\mathcal L}({\mathfrak g})$. We will denote the elements of this subgroup
by $p_\eta(x)$ for $x\in\bc$. 
\end{dfn}

The parametrizations can be chosen to have the following property:
for any $w\in\Waff$, any (real) root $\eta$ and any complex number 
$a$, $w p_{\eta}(a)w^{-1} = p_{w(\eta)}(a)$. 
For any two subgroups $\hat\cp\supset\hat\cq\supset\hat\cb$ such that $\Waff_{\hat\cp}$ is finite, let us 
consider the affine open neighborhood of $1$ in $\hat\cp/\hat\cq $,
$$
\cu^-(\hat\cp/\hat\cq) = \prod_{\substack{\eta <0, \eta {\rm\ real} \\ \cu_\eta\subset\hat\cp,\  \cu_\eta\not\subset \hat\cq}} \cu_\eta.
$$ 
If $w \in \Waff_{\hat\cp}/\Waff_{\hat\cq}$ (again we identify $w$ with a representative of minimal length in its class), 
let us consider the following sets of real roots 
$$
R^+(w) = \set{ \eta>0\mid \cu_{w^{-1}(\eta)}\not\subset \hat\cq}\quad {\rm and}\quad
R^-(w) = \set{ \theta<0\mid w(\theta) < 0,\ \cu_\theta\in \cu^-(\hat\cp/\hat\cq)},
$$
then $w \cu^-(\hat\cp/\hat\cq) = \cu^+(w) w \cu^-(w)$ is an affine open neighborhood of the point $w\in\hat\cp/\hat\cq$, where 
$$
\begin{array}{cc}
\cu^+(w) = \prod_{\eta\in R^+(w)} \cu_\eta  & \hbox{ and }\ 
\cu^-(w) = \prod_{\theta \in R^-(w)} \cu_\theta. 
\end{array}
$$ 
Note that if $w'$ is any element in the class of $w$ in $\Waff_{\hat\cp}/\Waff_{\hat\cq}$, then $R^+(w) = R^+(w')$. In addition, the set of positive real roots $R^+(w)\amalg \big (-R^-(w)\big )$ indexes the walls of $\tta$ that contain $F_\cp$, but not $w F_\cq$ (meaning that the set of associated ``affine'' roots, in $\tta$, does).
\begin{dfn}\rm
Let $\de\in\Gamma(\g_\lam)$ be a combinatorial gallery of type $t_{\gamma_\lam}$:
$$
\de =
[\de_0,\de_1,...,\de_r] = (F_f = \De'_0\subset \Delta_0
\supset\cdots\Delta_{j-1}\supset\Delta'_{j}\subset\Delta_{j}
\supset\cdots  \subset \De_p \supset F_\nu).
$$ 
We define a subset $\cu_\de$ of $\BSlam$ as follows, 
$$
\cu_\de  = \set{[g_0,g_1,...,g_p]\in\BSlam\mid\ g_j\in \cu^+(\de_j)\de_j \cu^-(\de_j)}.
$$ 
\end{dfn}
Now $\BSlam$ is a sequence of locally Zariski-trivial fibrations,
having partial flag varieties as fibres. So $\cu_\de$ is actually 
an {\it affine open subset of $\BSlam$ centered at $\de$}, and 
$\cu_\de$ is isomorphic to $\cu^+(\de_0)\de_0 \cu^-(\de_0)\times
\cdots\times \cu^+(\de_p)\de_p \cu^-(\de_p)$. 

\begin{rem}\rm
For any combinatorial gallery $\de$, the cell $C(\de)=\hat{r}_{\gamma_\lambda}^{-1}(\de)$ 
(Proposition~\ref{BBcell}) is a subvariety of the open set $\cu_\de$.
\end{rem}


\section{Minimal Galleries}\label{minimalgalleries}
Recall that the definition of a minimal gallery is actually valid in any apartment. 
Further, given any two faces of the building, there exists at least one apartment 
that contains those two faces, and any minimal gallery between them is contained 
in this apartment. Recall also that the Kac-Moody group acts on the affine building via simplicial maps.

To characterize the minimal galleries as a subset of $\BSlam$, let us begin with the 
case of a minimal gallery of alcoves. For the following proposition see \cite{BT}, 2.1.9: 

\begin{prop}\label{minimalalcovesgallery}
Let $\mu = (\De_f,\De_1,...,\De_t)$ be a sequence of alcoves in $\cj^{\mathfrak a}$. 
Then $\mu$ is a minimal gallery of alcoves if and only if there exists $b_j\in\hat\cb$ and 
$s_j\in \Saff$, for $j=1,...,t$, such that
$$
\De_j = b_1s_1\cdots b_js_j \De_f
$$ 
and $s_1s_2\cdots s_t$ is a reduced decomposition of an element in $\Waff$.
\end{prop}

\begin{rems}\rm 
1) By an abuse of notation we denote by $s_j$ the simple reflection in $\Waff$ as well as the 
representative $s_j = p_{\al_j}(1)p_{-\al_j}(-1)p_{\al_j}(1)$ in the Kac-Moody group (here $\al_j$ is the simple root in the root system of $\hat{\mathcal L}(G)$ corresponding to $s_j$).  

2) The $b_j$ are elements of the root subgroups, more precisely, 
$b_j = p_{\eta_j}(a_j)$ for some $a_j\in\bc$. 
Note that any element in $\cu^+(s_1\cdots s_j)$ can be written 
as $p_{\eta_1}(a_1)p_{s_1(\eta_2)}(a_2)\cdots p_{s_1\cdots s_{j-1}(\eta_j)}(a_j)$ 
for some complex numbers $a_i$'s. 

3) Minimal galleries of alcoves not starting at $\De_f$ can be obtained using the action of the Kac-Moody group.

4) In a minimal gallery of alcoves (of length $t$), any two consecutive alcoves are separated by a wall, and all 
these walls are distinct, so the gallery crosses exactly $t$ distinct walls in any apartment containing it.
\end{rems}

We come back to the general case. We cut the gallery of types into smaller blocks,    
i.e., we consider the following gallery of types $\varepsilon = (\theta\subset \kappa'\supset \kappa) = 
(t_{j-1}\subset t'_j\supset t_j)$ for some $j\in\set{1,...,p}$, where 
$t_{\g_\lam}= (S = t'_0\supset t_0 \subset \cdots\supset
t_{j-1} \subset t'_j \supset t_j \subset\cdots \supset t_p \subset t_{\lam})$. 
Let $\hat\cp$ and $\hat\cq$ be the parabolic subgroups containing $\hat\cb$ of type $\kappa'$ and $\kappa$, respectively. Let also $\cp$ and $\cq$, the corresponding parahoric subgroups. 
We know that $\hat\cp/\hat\cq = \bigcup_{w\in \Waff_{\hat\cp}/\Waff_{\hat\cq}} \cu^+(w)w\cu^-(w)$. We denote by $\tau^m$ the smallest 
representative of the class of the longest element in $\Waff_{\hat\cp}/\Waff_{\hat\cq}$. In addition, as a corollary of the Proposition~\ref{minimalalcovesgallery} and the remarks above, we get:

\begin{lem}
Let $g = (E\supset F'\subset F)$ be a gallery of type $\varepsilon$, starting at a 
face $E$ of the fundamental alcove $\De_f$. If $F = xF_\cq$, with $x\in\cu^+(\tau^m)\tau^m$, then $g$ is a minimal gallery.
\end{lem}

\proof 
By Proposition~\ref{minimalalcovesgallery}, we can construct a
minimal gallery of alcoves of length $r = \ell(\tau^m)$ between $\De_f\supset E$ and the alcove $x\De_f = proj_F(\Delta_f)$.
This shows that $\Delta_f$ is at maximal distance from $F$. Since any such minimal gallery is obtained from one of the form 
above using the operation of the stabilizer of $E\cup F$, the implication follows (see Lemma \ref{Stab} and Remark \ref{equivminisgall}).
\endpf

We will also need the following partial converse statement.

\begin{prop}\label{minigallery}
If a gallery $g = (E\supset F'\subset F)$ of type $\varepsilon$, starting at a 
face $E$ of the fundamental alcove $\De_f$, is minimal then 
$F = y'xF_\cq$ where $x\in \cu^+(\tau^m)\tau^m$ and $y'\in \hat\cq_E = Stab (E)$ 
(a subgroup of $\hat\cp$).
\end{prop}

\proof
Let $A$ be an apartment containing $g$. 
By definition, the gallery is minimal in $A$ if, for any alcove $\De\supset E$ at maximal distance to 
$F$ in $A$, the walls separating $\De$ from $F$ are exactly the ones that contain $F'$ but not $F$. 
Let us denote by $\sh(F',F)$ this set of walls and by $r$ its cardinality. 

Denote $\De\supset E$ an alcove at maximal distance to $F$ in $A$. Since
$g = (E\supset F'\subset F)$ is minimal, any minimal gallery of 
alcoves between $\De\supset E$ and the alcove $proj_F(\Delta)$ 
is of length $r$. Let $\mu = (\De = \De_0,\De_1,...,\De_r = proj_F(\De))$ be 
such a gallery. Since $\De\supset E\subset \De_f$, there 
exists an element $y\in \hat\cq_E$ such that $y\mu = (\De_f = y\De_0,y\De_1,...,y\De_r)$ 
is a minimal alcoves' gallery starting at $\De_f$. Hence, by 
Proposition~\ref{minimalalcovesgallery}, there exists $x\in\cu^+(\tau^m)\tau^m$ 
such that $y\ proj_F(\De) = x\De_f$. In particular, $F = y^{-1}xF_\cq$, hence the proposition is proved.
\endpf

Let $w\in \Waff_{\hat\cp}/\Waff_{\hat\cq}$, by abuse of notation we still
denote by $w$ the unique minimal representative in $\Waff_{\hat\cp}$ of its class. 
The following proposition expresses morally the ``trace'' of minimal galleries on the 
open set $\cu^+(w)w\cu^-(w)$ of $\hat\cp/\hat\cq$. We use the results of Deodhar \cite{Deo}.

\begin{prop}\label{traceofminimalgalleries}
The intersection
$\cu^+(w)w\cu^-(w) \cap\cu^+(\tau^m)\tau^m$ decomposes into a finite disjoint union 
of subvarieties, each of which is isomorphic to a product of $\bc$ and $\bc^*$.
\end{prop}

\proof
Let $w_\cp$ be the longest element in $\Waff_{\hat\cp}$, then $\cu^+(\tau^m)\tau^m = 
\cu^+(w_\cp)w_\cp$ modulo $\cq$. Next let $v\in\Waff_{\hat\cp}$ be such that $w_\cp = wv$ 
and $\ell(w) + \ell(v) = \ell(w_\cp)$, then $\cu^+(w)\cu^+(w_\cp)w_\cp = \cu^+(w_\cp)w_\cp$. 
Hence we have:
$$
\cu^+(w)w\cu^-(w) \cap\cu^+(\tau^m)\tau^m = \cu^+(w)\big ( w\cu^-(w)\cap \cu^+(w_\cp)w_\cp\big ).
$$
Let $\hat\cb^-$ is the subgroup of the Kac-Moody group $\hat{\mathcal L}(G)$ 
generated by the torus $\ct$ and all the root subgroups associated to negative 
roots, then   
$$\begin{array}{rcl}
w\cu^-(w)\cap \cu^+(w_\cp)w_\cp & = & w\cu^-(w)\cap w_\cp\cu^-(\hat\cp/\hat\cq)\\
 & = & w_\cp\big (w_\cp^{-1}w\cu^-(w)\cap \cu^-(\hat\cp/\hat\cq)\big )\\
 & = & w_\cp\big (\cu^+(v^{-1})v^{-1}\cap \cu^-(\hat\cp/\hat\cq)\big )\\
 & = & w_\cp\big (\hat\cb v^{-1}\cap \hat\cb^-\hat\cq/\hat\cq\big ),\\
\end{array}$$
By \cite{Deo}, the intersection $\hat\cb v^{-1}\cap \hat\cb^-\hat\cq/\hat\cq$ in $\hat\cp/\hat\cq$ 
decomposes into a finite disjoint union of subvarieties isomorphic to products of 
$\bc$ and $\bc^*$, which implies  the proposition.
\endpf

\begin{rem}\label{CstarOpen}
\rm
In the Deodhar decomposition of $\hat\cb v^{-1}\cap \hat\cb^-\hat\cq/\hat\cq$, there is exactly one 
subvariety of maximal dimension which consists only of copies of $\bc^*$ 
(this is the one which corresponds to the subexpression $(1,1,...,1)$ of $v^{-1}$). In our situation, 
this is exactly the open and dense subset  where the variables associated to the root subgroups 
in $\cu^-(w)$ are taken to be nonzero:
$$
\cu^+(w)w\cu^-(w)(\bc^*) \subsetneq \cu^+(w)w\cu^-(w) \cap\cu^+(\tau^m)\tau^m.
$$ 
\end{rem}


\section{The Cell Associated to a Positively Folded Gallery}\label{cell}

In this section, we describe the fibre $C(\de) = \hat r_{\g_\lam}^{-1}(\de)$ of
the retraction $\hat r_{\g_\lam} : \BSlam\to \G(\g_\lam)$ over a positively folded 
gallery $\de\in\G^+(\g_\lam,\nu)$, as a subvariety 
of the affine open subset $\cu_\de$ of $\BSlam$. Let us fix  
$$
\de = [\de_0,\de_1,...,\de_p] = (F_f = \De'_0\subset \Delta_0\supset\Delta'_1\subset\cdots\Delta_{j-1}\supset\Delta'_{j}\subset\Delta_{j}\supset\cdots  \subset \De_p \supset F_\nu).
$$  
For $j\in\set{0,...,p}$, we denote by $\tau^m_j$ the ``element of maximal length'' 
in $W'_j/W_j$, let $\ell_j$ be the length (see Example \ref{taum} in Section \ref{combinatorial gallery}). We fix a numbering on the set of walls containing $F_{\cp_j}$, 
but not $\de_j F_{\cq_j}$ (note that this set is equal to $(\de_0\cdots\de_{j-1})^{-1}\sh_j$). The numbering
is chosen such that the indexing set can be decomposed into
$$
I_j = I^+_j\amalg I^-_j = \set{j_1,...,j_{\ell(\de_j)}}\amalg\set{j_{\ell(\de_j)+1},...,j_{\ell_j}},
$$ 
so that $R^+(\de_j) = \set{\eta_{j_i}\mid i\in I^+_j}$ and $R^-(\de_j) = \set{\theta_{j_i} \mid i\in I^-_j}$
(notation as in Section~\ref{opencovering}). 

Let $J_{-\infty}(\de)\subset I_0\amalg\cdots\amalg I_p$ be the subset of all the 
indices corresponding to a load-bearing wall of $\de$. We have a decomposition 
$J_{-\infty}(\delta) = J_{-\infty}^+(\delta)\amalg J_{-\infty}^-(\delta)$, 
where we define 
$$
J_{-\infty}^+(\delta) 
= J_{-\infty}(\delta)\cap\big (I_0^+\amalg\cdots\amalg I^+_p\big ) 
= \bigcup_{j=0}^p \big(J_{-\infty}(\delta)\cap I_j^+) ,
$$  
and similarly for $J_{-\infty}^-(\delta)$. Observe that $J_{-\infty}^-(\delta) = I_0^-\amalg\cdots\amalg I^-_p$
since the gallery $\de$ is positively folded, and $\dim\de = \sharp J_{-\infty}(\delta)$ 
(see Section \ref{dimensiongallery} for the definition of $\dim\de$). 

\begin{lem}\label{InvRetractNonreg}
Let $\de = [\de_0,\de_1,...,\de_p]\in\Gamma^+(\g_\lam)$ be a positively folded gallery.
Then $C(\de)= \hat r_{\g_\lam}^{-1}(\de)\subset \BSlam$ is an affine cell of dimension
$\dim\de$. It is the subvariety of $\cu_\de$ consisting of all the galleries
$[g_0,g_1,...,g_p]$ such that
$$
g_j \in \Big (\prod_{\ i\in J_{-\infty}(\de)\cap I^+_j} \cu_{\eta_{j_i}}\Big )\de_j \prod_{\ i\in I_j^-} \cu_{\theta_{j_i}}.
$$
\end{lem}

The fibre $r_{\gamma_\lam}^{-1}(\delta)\subset\cg_\lam$ is obtained as the intersection
$ \hat r_{\gamma_\lam}^{-1}(\delta)\cap \hat\cg_\lam$, so we have:
\begin{coro}\label{fibrequasiaffine}
The fibre over a positively folded gallery $\de$ of the map
$r_{\gamma_\lam}:\cg_\lam\to\Gamma^+(\gamma_{\lam})$
is naturally endowed with the structure of an irreducible quasi-affine variety. More precisely,
it is an open dense subset of an affine space of dimension $\dim\de$.
\end{coro}
Before we proceed with the proof of the lemma, we shall discuss the case where $\lam$ is a {\it regular}
dominant co-character. In this case, $\de = [\de_0,\de_1,..., \de_p]$ is such that
$\de_0\in W$ and $\de_j = 1$ or $s_{i_j}$, a simple reflection associated to the simple
(Kac-Moody) root $\al_{i_j}$. Equivalently, the $\De_j$'s are alcoves and the $\De'_j$'s are faces of
codimension 1 in two consecutive alcoves. Further, the set $I_j$ is reduced to a single element, it
corresponds to the wall which is the support of $\De'_j$. So, in case $\lam$ is regular, the lemma
has the following form.

\begin{lem}
Assume $\lam$ regular.
Let $\de = [\de_0,\de_1,...,\de_p]\in\Gamma^+(\g_\lam)$ be positively folded gallery.
Then $C(\de)$ is the subvariety of $\cu_\de$ consisting of all the galleries
$[g_0,g_1,...,g_p]$ such that
$$
g_0 \in \de_0\prod_{\be <0, \de_0(\be) < 0} \cu_\be\hbox{ and }
g_j = \left\{
\begin{array}{l}
\delta_j \hbox{ if } j\not\in J_{-\infty}(\delta)\\
p_{-\alpha_{i_j}}(a_j),\,a_j\in\bc, \hbox{ if } j\in J_{-\infty}^-(\delta)\\
p_{\alpha_{i_j}}(a_j)s_{i_j},\,a_j\in\bc, \hbox{ if } j\in J_{-\infty}^+(\delta).\\
\end{array}
\right.
$$
\end{lem}
\proof
Consider first the condition on $g_0$. We want to characterize the alcoves
of the form $g_0(\De_f)$ containing $0$ which retract onto $\de_0(\De_f)$
by $\rinf$. In this case we can identify the retraction
centered at $-\infty$ with the retraction onto the apartment $\Aaff$ centered at
$w_0(\De_f)$, i.e., $\rinf\big (g_0(\De_f)\big ) = r_{w_0(\De_f),\Aaff}\big (g_0(\De_f)\big )$.
So the problem amounts to find all the spherical chambers that retract on the
chamber $\de_0(\Lc_f)$ via the retraction $r_{-\Lc_f} = r_{w_0(\Lc_f),\A^s}$ onto the
spherical apartment, centered at $w_0(\Lc_f) = -\Lc_f$. But we know that
$r_{-\Lc_f}^{-1}\big (\de_0(\Lc_f)\big ) = r_{-\Lc_f}^{-1}\big (\de_0w_0(-\Lc_f)\big ) =
B^-\de_0w_0$ in $G/B^-$ (recall Example \ref{exampleBruhat}). In addition,
$$
B^-\de_0w_0 = \Big (\prod_{\eta <0, (\de_0w_0)^{-1}(\eta) > 0} \cu_\eta\Big )\de_0w_0
= \de_0\prod_{\be <0, \de_0(\be) < 0} \cu_\be \ w_0.
$$
Using the isomorphism $G/B^-\simeq G/B$ given by right multiplication by $w_0$,
we see that the alcoves retracting onto $\de_0(\De_f)$ via $\rinf$ are of the form
$g_0(\De_f)$, $g_0\in \de_0\prod_{\be <0, \de_0(\be) < 0} \cu_\be$.

The conditions on $g_j$ for $j\in\set{1,...,p}$ follow easily from the definition of the
retraction $\rinf$. Let $g = [g_0,...,g_p]$ be a gallery, we want to know when $g$
retracts onto $\de$. We can assume that $g$ has been retracted
until the index $j-1$, i.e., 
$$
g = [\de_0,\de_1,...,\de_{j-1},g_j,...,g_p] = (F_f\subset \Delta_0\supset\Delta'_1
\subset\cdots\Delta_{j-1}\supset \Delta'_{j} \subset F_{j}\supset\cdots
\supset F'_\nu),
$$
with $F_{j} = \de_0\de_1\cdots\de_{j-1} g_{j}(\De_{f})$. If $F_j=\Delta_j$, then this corresponds to $a_j=0$ in the first or the second condition. Otherwise the gallery
$(\Delta_{j-1}\supset \Delta'_{j} \subset F_{j})$ (resp. $(s_{\ch_j}\Delta_{j-1}\supset \Delta'_{j} \subset F_{j})$) is minimal in any apartment containing it, where $s_{\ch_j}$ is the reflection defined by the wall $\ch_j$ supporting the face $\De'_j$.
Let $A$ be an apartment containing $F_j$, the wall $\ch_j$ and a representative sector in the equivalence
class of $\mathfrak C_{-\infty}$ in $\tta$.
Let us incarnate the retraction center at $-\infty$ into $\rinf = r_{\Sigma,\tta}$, where
$\Sigma\in\tta\cap A$ is a sufficiently far away alcove from $\De'_j$.

If $j\not\in J_{-\infty}(\de)$, then $\ch_j$ is not load-bearing at this place. In particular, $\Sigma$ and $\De_j$ are not
separated by $\ch_j$. But $\Sigma$ and $F_j$ are separated by $\ch_j$ in $A$.
Hence, because of the properties of the retraction, $r_{\Sigma,\tta}(F_j) \ne \De_j$, unless
$g_j = \de_j$ (and $F_j = \De_j$). Now, if $j\in J_{-\infty}^-(\de)$, the wall $\ch_j$ is
load-bearing and $\de_j = 1$. So this time the only condition on $g_j\in\hat\cp_j/\hat\cq_j
= \cu^-(\de_j)\hat\cq_j/\hat\cq_j\amalg s_{i_j}\hat\cq_j/\hat\cq_j$ is that $g_j\ne s_{i_j}$. Hence,
$F_j$ retracts on $\De_j$ if and only if $g_j = p_{-\al_{i_j}}(a_j)$, for some complex
number $a_j$. The case where $j\in J_{-\infty}^+ (\de)$ is treated similarly and
corresponds to the decomposition $\hat\cp_j/\hat\cq_j = \cu^+(\de_j) s_{i_j}\hat\cq_j/\hat\cq_j\amalg \hat\cq_j/\hat\cq_j$,
this finishes the proof of the lemma.
\endpf

{\it Proof of Lemma \ref{InvRetractNonreg}.}
The proof for $g_0$ is similar to the proof in the regular case,
the computation takes now place in the quotients $G(\co)/\cq_0 \simeq G/P_\lam$ and $G/P_\lam^-$,
where $P^-_\lam$ is the parabolic subgroup of type $t_0$ containing $B^-$. So, applying
the retraction $r_{-\Lc_f} = r_{w_0(\Lc_f),\A^s}$ to the spherical faces of type $t_0$, we obtain that $g_0F_{\cq_0}$
retracts onto $\de_0F_{\cq_0}$ if and only if
$$
g_0\in \de_0\prod \cu_\be,\ \hbox{\rm where\ }
\beta\in R^-(\de_0)=\{\be <0,\ \de_0(\be) < 0,\, \cu_\be\in \cu^-(G/P_\lam)\}.
$$
Note that  $J_{-\infty}(\de)\cap I_0^+ = \emptyset$.
Next consider the case $j>0$. Let $g = [g_0,g_1,...,g_p]\in\BSlam$ be a gallery.
Again, we can retract step by step, i.e. we can assume that $g$ has been retracted
to a gallery $g'$ up to  the index $j-1$, so
$$
g' = [\de_0,\de_1,...,\de_{j-1},g_j,...,g_p] = (F_f\subset \Delta_0\supset\Delta'_1\subset
\cdots\Delta_{j-1}\supset \Delta'_{j} \subset F_{j}\supset F'_{j+1}\cdots  \supset F'_\nu),
$$
with $F_{j} = \de_0\de_1\cdots\de_{j-1} g_{j} F_{\cq_j}$ and $g_j\in \cu^+(\de_j)\de_j \cu^-(\de_j)$.
Hence, there exists some complex numbers $a_i$ for $i\in I^+_j$, and $b_k$ for $k\in I^-_j$ such
that we can write $g_j$ as follows
\begin{equation}\label{writingofg_j}
g_j = p_{\zeta_1}(a_1)s_{\zeta_1}\cdots p_{\zeta_{\ell(\de_j)}}
(a_{\ell(\de_j)})s_{\zeta_{\ell(\de_j)}}\prod_{k\in I^-_j} p_{\theta_{j_k}}(b_k),
\end{equation}
where $s_{\zeta_1}\cdots s_{\zeta_{\ell(\de_j)}}$ is a reduced
decomposition of $\de_j$ and
$$
R^+(\de_j) = \set{\zeta_1, s_{\zeta_1}(\zeta_2),..., s_{\zeta_1}
\cdots s_{\zeta_{\ell(\de_j)-1}}(\zeta_{\ell(\de_j)}) }.
$$
Assume first $b_k=0$ for all $k$. By
Proposition \ref{minimalalcovesgallery}, the alcove gallery
$\mu = (\G_0,\G_1,\G_2,...,\G_{\ell(\de_j)})$,
where $\G_i = \de_0\de_1\cdots\de_{j-1} p_{\zeta_1}(a_1)s_{\zeta_1}
\cdots p_{\zeta_{i}}(a_{i})s_{\zeta_{i}} \De_f$, is a minimal gallery between
$\G_0\supset\De_{j-1}$ and $\G_{\ell(\de_j)}\supset F_j$. Because of the
minimality, $F_j$ retracts on $\De_j$ if and only if $\mu$ retracts onto the minimal
gallery $\xi = (\Xi_0 , \Xi_1 ,..., \Xi_{\ell(\de_j)})$ in $\tta$, where
$\Xi_i = \de_0\de_1\cdots\de_{j-1} s_{\zeta_1}\cdots s_{\zeta_{i}}\De_f$. Again, we can
assume that $\mu$ has been retracted onto $\xi$ until the index $i-1$, so we have to
find out how the alcove $\G'_i = \de_0\de_1\cdots\de_{j-1} s_{\zeta_1}
\cdots s_{\zeta_{i-1}}p_{\zeta_{i}}(a_{i})s_{\zeta_{i}} \De_f$ may retract onto $\Xi_i$. But we are in the situation where we can apply the discussion concluding the previous proof. Hence, if $i\in J_{-\infty}(\de)\cap I^+_j$, $a_i$ can be any complex number, otherwise $a_i$ must equal $0$. Therefore, we get the
expected conditions on the variables associated to indices in $I^+_j$.

To conclude, the fact that there is no condition on the $b_k$'s in (\ref{writingofg_j}) is a
consequence of an analogue computation as in the case $j=0$ and
the assumption that the gallery is positively folded. Thus, the lemma is proved.
\endpf

\begin{rem}\rm
The definition of load-bearing wall with respect to $-\infty$ makes sense for any gallery 
and the previous description of the fibre also extends to arbitrary combinatorial 
galleries. Therefore, one gets a description in terms of galleries of the 
Bia{\l}ynicki-Birula cellular decomposition of the Bott-Samelson variety $\BSlam$.
\end{rem}


\section{LS-galleries and MV-cycles}\label{mvcycles}

We describe an affine open subset of $r_{\g_\lam}^{-1}(\de)=C(\delta)\cap \cg_\lam$, 
isomorphic to a product of $\bc$'s and $\bc^*$'s.  Thereby, we can associate to each LS-gallery 
$\de$ an open dense subset of a Mirkovic-Vilonen cycle $Z(\de)$.  We first deal 
for simplicity with the {\it regular} case.

\begin{prop}
Let $\lam$ be a regular dominant co-character and
let $\de = [\de_0,\de_1,...,\de_p]\in\Gamma^+(\g_\lam)$. 
Then $r_{\g_\lam}^{-1}(\de)$ is a subvariety of $\cu_\de$ isomorphic to a product of
$\bc$'s and $\bc^*$'s. More precisely, the fibre consists of all the galleries
$[g_0,g_1,...,g_p]$ such that
$$
g_0 \in \de_0\prod_{\be <0, \de_0(\be) < 0} \cu_\be\hbox{ and }
g_j = \left\{
\begin{array}{l}
\delta_j \hbox{ if } j\not\in J_{-\infty}(\delta)\\
p_{-\alpha_{i_j}}(a_j), \ a_j\ne 0 \hbox{ if } j\in J_{-\infty}^-(\delta),\\
p_{\alpha_{i_j}}(a_j)s_{i_j} \hbox{ if } j\in J_{-\infty}^+(\delta).\\
\end{array}
\right.
$$
\end{prop}

\proof
Combining the Propositions~\ref{minimalalcovesgallery} and \ref{minigallery}, we
get the following simple criterion in the regular case: 
a gallery of alcoves $(\De_f\supset F'\subset \De)$, where $F'$ is a face of
codimension~1 associated to the simple root $\al$, is minimal if and only if $\De = x\De_f$ 
with $x \ne 1$ in $\cp_\al/\cb$. This applies to each step in the gallery $[g_0,g_1,...,g_p]$, 
and the only condition to add in order for the latter to be minimal is the one of the lemma.
\endpf

In the general case the situation is slightly more complicated, but we still have:

\begin{thm}\label{mvcycledescription}
Let $\lam$ be a dominant co-character and let $\nu\prec\lam$ be an arbitrary co-character.
For each positively folded gallery $\de\in \G^+(\g_\lam,\nu)$, 
the fibre $r_{\g_\lam}^{-1}(\de)$ admits a decomposition into a union of subvarieties,
each of which is isomorphic to a product of $\bc$'s and $\bc^*$'s. 
In particular, the subvariety $O(\de)\simeq \bc^{a}\times(\bc^*)^{b}$ of $r_{\g_\lam}^{-1}(\de)$ 
is open and dense, where $a+b=\dim\de$, $b=\sharp J_{-\infty}^-(\delta)$, and
$O(\de)$ consists of all the galleries $[g_0,g_1,...,g_p]$ such that
$$
g_j \in \Big (\prod_{\ i\in J_{-\infty}(\de)\cap I^+_j} \cu_{\eta_{j_i}}\Big )\de_j 
\prod_{\ k\in I_j^-} \cu_{\theta_{j_k}}(\bc^*).
$$ 
\end{thm}

By \cite{MV},  we know that $U^-(\ck)\nu\cap \cg_\lam$ is 
equidimensional. It follows that $\overline{U^-(\ck)\nu\cap \cg_\lam}$ is the
union of the $Z(\delta)=\overline{r_{\g_\lam}^{-1}(\de)}$ for $\de\in \G^+(\g_\lam,\nu)$
of maximal dimension, i.e., $\de$ a LS--gallery.

\begin{coro} The MV-cycles corresponding to $\nu$ 
(i.e., the irreducible components of $\overline{U^-(\ck)\nu\cap \cg_\lam}$), 
are given by the closures $Z(\de)=\overline{O(\de)}$ for $\de\in\Gamma_{LS}^+(\g_\lam,\nu)$.
\end{coro}

{\it Proof of the theorem.\/}
The first definition of $\BSlam$ (Definition \ref{defiBott1}) allows to iterate the use 
of the Propositions \ref{minigallery} and \ref{traceofminimalgalleries}. 
The theorem is then obtained using Remark \ref{CstarOpen} and the same arguments
as above in the regular case.
\endpf

It is not yet clear how to express the inclusion relations between the closures
of the fibres in terms of the combinatorics of the galleries.
For simplicity assume $\lam$ is regular and fix $w\in W$, then one sees easily for 
$\de'= [w',s_{i_1},...,s_{i_p}]$ and $\de = [w,s_{i_1},...,s_{i_p}]$:
$$
r^{-1}_{\g_\lam}(\de')\subset \overline{r^{-1}_{\g_\lam}(\de)}\Leftrightarrow w\le w'.
$$
More generally, let $\de,\g$ be two consecutive galleries in 
the sequence of positive foldings starting with $w\g_\lam = [w,s_{i_1},...,s_{i_p}]$, 
where the $s_{i_j}$'s are the simple reflections associated to the simple (Kac-Moody) 
root $\al_{i_j}$ (for the sequence of foldings, see Section \ref{combinatorial gallery} (\ref{sequenceoffoldings})). More precisely, there exists an index $k$ such that
\begin{equation}\label{foldingclosure}
\de =  [w,\de_1,...,\de_{k-1},  s_{i_k},  s_{i_{k+1}},..., s_{i_p}]\quad \hbox{\rm and}\quad
\g = [w,\de_1,...,\de_{k-1},   1 ,  s_{i_{k+1}},..., s_{i_p}].
\end{equation}

\begin{prop}\label{fibreclosure}
If $\de,\g$ are as in (\ref{foldingclosure}), then $r^{-1}_{\g_\lam}(\de) \subset \overline{r^{-1}_{\g_\lam}(\g)}$.
\end{prop}

\proof
Let $g = [g_0,g_1,...,g_{k-1}, p_{-\al_{i_k}}(a_k), g_{k+1},...,g_p]$ be a 
point in $r^{-1}_{\g_\lam}(\g)$, where $g_j = p_{\al_{i_j}}(a_j) s_{i_j}$, 
for $j\geq i+1$ and some complex number $a_j$ (which might be zero). 
For any positive (Kac-Moody) root $\al$ and any complex number $a$, one has
$p_{-\al}(a) = p_{\al}(a^{-1})s_\al h_\al p_{\al}(a^{-1})$ (see \cite{T87}) for some torus element $h_\al$. 
So, $g$ may be written as
$$
g = [g_0,g_1,...,g_{k-1}, p_{\al_{i_k}}(a_k^{-1})s_{i_k} h_{\al_{i_k}} p_{\al_{i_k}}(a_k^{-1}), g_{k+1},...,g_p].
$$ 
But the last part of this gallery is minimal, so it is stable by $p_{\al_{i_k}}(a_k^{-1})$, and hence
$$
g = [g_0,g_1,...,g_{k-1}, p_{\al_{i_k}}(a_k^{-1})s_{i_k}, g'_{k+1},...,g'_p],
\quad\hbox{\rm where\ }g'_j = p_{\al_{i_j}}(a'_j) s_{i_j},\,a'_j = b_j a_j
$$ 
for $j\geq i+1$  
and some $b_j\in\bc^*$. Therefore, $r^{-1}_{\g_\lam}(\de) \subset \overline{r^{-1}_{\g_\lam}(\g)}$.
\endpf

\begin{coro}\label{corofibreclosure}
If $\de,\g$ are two consecutive LS--galleries in a sequence of foldings, 
then the associated MV--cycles are contained in each other: $Z(\de)\subset Z(\g)$.
\end{coro}

\begin{rem}\rm
The same kind of arguments (using the open subvariety $O(\delta)\subset r^{-1}_{\g_\lam}(\de)$) apply 
also in the general case, i.e., Proposition~\ref{fibreclosure} and Corollary~\ref{corofibreclosure}
hold for an arbitrary dominant co-character $\lam$ and $\de\in\G^+(\g_\lam)$. 
\end{rem}
\begin{exam}
\rm
Let us take a simple example to illustrate this section. Recall the setting of Example \ref{exampleSl3}: 
$G = SL_3(\bc)$, $\Saff = \set{s_0,s_1,s_2}$ and we let $\lam = \be = \al_1 +\al_2$ be the highest root. 
We set $\g_\be = [1,s_0]$. The non-minimal positively folded galleries in $\G^+(\g_\be)$ are 
$$
\de_{12} = [s_1s_2, 1],\ \de_{21} = [s_2s_1, 1], \de_{121} = [s_1s_2s_1, 1].
$$ If we apply the theorem, we get:
$$
r_{\g_\be}^{-1}(\de_{12}) = \set{[s_1s_2p_{-\al_1}(a), p_{-\al_0}(b)]\in\hat\Sigma(t_{\g_\be}),\ b\ne 0},
$$
$$
r_{\g_\be}^{-1}(\de_{21}) = \set{[s_2s_1p_{-\al_2}(a), p_{-\al_0}(b)]\in\hat\Sigma(t_{\g_\be}),\ b\ne 0},
$$ and
$$
r_{\g_\be}^{-1}(\de_{121}) = \set{[s_1s_2s_1, p_{-\al_0}(b)]\in\hat\Sigma(t_{\g_\be}),\ b\ne 0}.
$$ 
The two Mirkovic-Vilonen cycles associated to the pair $(\lam = \be,\nu = 0$) are obtained as the 
closure of the two varieties $r_{\g_\be}^{-1}(\de_{12})$ and $r_{\g_\be}^{-1}(\de_{21}) $, and
one can show that $r_{\g_\be}^{-1}(\de_{121})$ is in the closure of those two.
Using Proposition~\ref{fibreclosure} and, for the inclusions of the $Z(\de_{ij})$, some arguments similar to 
the ones above, we get the following inclusion relations (note the very suggestive similarity of the diagram on
the left with the crystal graph):
$$
\begin{array}{cccccc}
Z([s_1s_2, s_0])&\subset&Z(\de_{12})&\subset&Z([s_1, s_0])\\
\cup& & & &\cap\\
Z([w_0,s_0])&& & & Z([1, s_0])\\
\cap& & & & \cup\\
Z([s_2s_1, s_0])&\subset&Z(\de_{21})&\subset&Z([s_2, s_0])\\
\end{array}
\quad\hbox{\rm and}\quad
\begin{array}{c}
Z([s_2s_1, s_0])\subset Z([s_1, s_0])\\
\\
Z([s_1s_2, s_0])\subset Z([s_2, s_0]).\\
\end{array}
$$
\end{exam}


\end{document}